\documentclass{gtart_h}  

\def\ifplaintex{\expandafter\ifx\csname documentclass\endcsname\relax}

\ifplaintex 
\else
\fi

\expandafter\ifx\csname beginpicture\endcsname\relax
\expandafter\ifx\csname documentclass\endcsname\relax
\input pictex \else
\input prepictex \input pictex \input postpictex \fi\fi

\def\gt{{\mathsurround=0pt\it $\cal G\mskip-2mu$eometry \&\ 
$\cal T\!\!$opology}}        

\def\gtp{{\mathsurround=0pt\it $\cal G\mskip-2mu$eometry \&\ 
$\cal T\!\!$opology $\cal P\!$ublications}}  

\def\lognumber#1{\def\thelognumber{#1}}
\def\volumenumber#1{\def\thevolumenumber{#1}}
\def\papernumber#1{\def\thepapernumber{#1}}
\def\volumeyear#1{\def\thevolumeyear{#1}}

\def\pagenumbers#1#2{\def\startpage{#1}\def\finishpage{#2}}
\def\published#1{\def\publishdate{#1}}
\def\proposed#1{\def\theproposer{#1}}
\def\seconded#1{\def\theseconders{#1}}
\def\received#1{\def\receiveddate{#1}}
\def\revised#1{\def\reviseddate{#1}}
\def\accepted#1{\def\accepteddate{#1}}

\def\asciiemail#1{\def\theasciiemail{#1}}

\long\def\asciiabstract#1{\long\def\theasciiabstract{#1}}

\let\\\par\let\thelognumber\relax
\let\thevolumenumber\relax\let\thepapernumber\relax
\let\thevolumeyear\relax\let\thesamplenumber\relax\let\startpage\relax
\let\finishpage\relax\let\publishdate\relax\let\receiveddate\relax
\let\reviseddate\relax\let\accepteddate\relax\let\theasciititle\relax
\let\theasciiauthors\relax
\let\theasciiabstract\relax
\let\theasciiemail\relax\let\theshortauthors\relax\let\theshorttitle\relax

\long\def\maketitlep{   

\count0=\startpage

\gt\hfill      
\beginpicture
\setcoordinatesystem units <0.33truein, 0.33truein> point at 2.2 0.9
\setplotsymbol ({$\cal G$})
\plotsymbolspacing=9truept
\circulararc 315 degrees from 0 1 center at 0 0
\setplotsymbol ({$\cal T$})
\circulararc 315 degrees from 1 -1 center at 1 0
\endpicture
\break
{\small\ifx\thesamplenumber\relax 
Volume \else Sample
\fi\thevolumenumber\ (\thevolumeyear)
\startpage--\finishpage\nl
Published: \publishdate}
\vglue 0.5truein plus 0.4fil minus 0.1truein

{\parskip=0pt\leftskip 0pt plus 1fil\def\\{\par\smallskip}{\ifplaintex\large
\else\Large\fi\bf\thetitle}\par\medskip}   

\vglue 0pt plus 0.1fil 

{\parskip=0pt\leftskip 0pt plus 1fil\def\\{\par}{\sc\theauthors}
\par\medskip}

\vglue 0pt plus 0.1fil 

{\small\parskip=0pt\let\newline\\
{\leftskip 0pt plus 1fil\def\\{\par}{\sl\theaddress}\par}
\expandafter\ifx\theemail\relax    
\relax\else\vglue 5pt plus 0.02fil minus 2pt\def\\{\stdspace{\rm 
and}\stdspace} 
\cl{Email:\stdspace\tt\theemail}\fi
\ifx\theurl\relax                  
\relax\else\vglue 5pt plus 0.02fil minus 2pt\def\\{\stdspace{\rm 
and}\stdspace}
\cl{URL:\stdspace\tt\theurl}\fi\par}

\vglue 7pt plus 0.3fil minus 3pt

{\bf Abstract}
\vglue 5pt plus 0.1fil minus 2pt

\theabstract

\vglue 7pt plus 0.3fil minus 3pt

{\bf AMS Classification numbers}\quad Primary:\quad \theprimaryclass

Secondary:\quad \thesecondaryclass

\vglue 5pt plus 0.3fil minus 2pt

{\bf Keywords:}\quad \thekeywords

\vglue 10pt plus 0.5fil minus 5pt

{\small  Proposed: \theproposer\hfill Received: \receiveddate\nl
Seconded: \theseconders\hfill 
\ifx\reviseddate\relax                         
Accepted: \accepteddate                        
\else
Revised: \reviseddate                          
\fi}
\eject
}       

\let\maketitlepage\maketitlep
\let\maketitle\maketitlepage

\font\phead=cmsl9 scaled 950
\font=cmsl9 scaled 1050
\font\pnum=cmbx10 scaled 913
\font\lnum=cmbx10 
\font\pfoot=cmsl9 scaled 950
\font=cmsl9 scaled 1050
\ifplaintex
\headline{\vbox to 0pt{\vskip -4.5mm\line{\small\phead\ifnum
\count0=\startpage ISSN 1364-0380 (on line)
1465-3060 (printed) \hfill {\pnum\folio}\else\ifodd\count0\def\\{ }%
\ifx\theshorttitle\relax\thetitle\else\theshorttitle\fi\hfill{\pnum\folio}
\else\def\\{ and }{\pnum\folio}\hfill\ifx\theshortauthors\relax\theauthors
\else\theshortauthors\fi\fi\fi}\vss}}
\footline{\vbox to 0pt{\vglue 0mm\line{\small\pfoot\ifnum\count0=\startpage
\copyright\ \gtp\hfill\else
\gt, Volume \thevolumenumber\ (\thevolumeyear)\hfill\fi}\vss
}}
\else
\makeatletter
\def\@oddhead{{\small\lhead\ifnum\count0=\startpage ISSN 1364-0380 (on line)
1465-3060 (printed) \hfill {\lnum\number\count0}\else\ifodd\count0
\def\\{ }\ifx\theshorttitle\relax \thetitle \else\theshorttitle\fi\hfill
{\lnum\number\count0}\else\def\\{ and }{\lnum\number\count0}
\hfill\ifx\theshortauthors\relax 
\theauthors\else\theshortauthors\fi\fi\fi}}\def\@evenhead{\@oddhead}
\def\@oddfoot{\small\lfoot\ifnum\count0=\startpage\copyright\ \gtp\hfill\else
\gt, Volume \thevolumenumber\ (\thevolumeyear)\hfill\fi}
\def\@evenfoot{\@oddfoot}
\makeatother
\fi

\newwrite\gtoutfile
\long\gdef\makeheadfile{  
{\def\\{, }\def\s{ }
\immediate\openout\gtoutfile head.xxx
\immediate\write\gtoutfile{Proxy-for: \ifx\theasciiauthors\relax
\theauthors\else\theasciiauthors\fi\s<\ifx\theasciiemail\relax\theemail\else\theasciiemail\fi>}
\immediate\write\gtoutfile{\noexpand\\}
\immediate\write\gtoutfile{Authors: \ifx\theasciiauthors\relax
\theauthors\else\theasciiauthors\fi}
{\def\\{ }\immediate\write\gtoutfile{Title: \ifx\theasciititle\relax
\thetitle\else\theasciititle\fi}}
\immediate\write\gtoutfile{Subj-class: GT or SG or MG etc}
\immediate\write\gtoutfile{MSC-class: \theprimaryclass\ifx\thesecondaryclass\relax\else, \thesecondaryclass\fi}
\immediate\write\gtoutfile{Journal-ref: Geom. Topol. \thevolumenumber
(\thevolumeyear) \startpage-\finishpage}
\immediate\write\gtoutfile{Comments: Published by Geometry and Topology at}
\immediate\write\gtoutfile{\s\s http://www.maths.warwick.ac.uk/gt/GTVol\thevolumenumber/paper\thepapernumber.abs.html}
\immediate\write\gtoutfile{\noexpand\\}
\immediate\write\gtoutfile{}
\ifx\theasciiabstract\relax
\immediate\write\gtoutfile{\theabstract}\else
\immediate\write\gtoutfile{\theasciiabstract}\fi
\immediate\write\gtoutfile{}
\immediate\write\gtoutfile{\noexpand\\}
\immediate\write\gtoutfile{}
\immediate\closeout\gtoutfile}}  

\def\maketitlepage{\maketitlep\makeheadfile}
\let\maketitle\maketitlepage

\lognumber{557}
\received{27 January 2005}

\volumenumber{9}\papernumber{49}\volumeyear{2005}
\pagenumbers{2159}{2191}   
\revised{11 November 2005}
\published{22 November 2005}
\accepted{16 October 2005}
\proposed{Robion Kirby}
\seconded{Cameron Gordon, Joan Birman}

\usepackage{amsmath,amscd,amssymb}

\def\a{\alpha}

\def\e{\epsilon}
\def\f{\frac}
\def\g{\gamma}
\def\G{\Gamma}

\def\lb\{{\left\{}
\def\la{\lambda}
\def\La{\Lambda}
\def\lla{\longleftarrow}
\def\lm{\limits}
\def\lra{\longrightarrow}
\def\dllra{\Longleftrightarrow}
\def\llra{\longleftrightarrow}
\def\n{\nabla}
\def\ngth{\negthickspace}
\def\ola{\overleftarrow}
\def\Om{\Omega}
\def\om{\omega}
\def\op{\oplus}
\def\oper{\operatorname}
\def\oplm{\operatornamewithlimits}
\def\ora{\overrightarrow}
\def\ov{\overline}
\def\ova{\overarrow}
\def\ox{\otimes}
\def\p{\partial}
\def\rb\}{\right\}}
\def\s{\sigma}
\def\sbq{\subseteq}
\def\spq{\supseteq}
\def\sqp{\sqsupset}
\def\supth{{\text{th}}}
\def\T{\Theta}
\def\th{\theta}
\def\tl{\tilde}
\def\thra{\twoheadrightarrow}
\def\un{\underline}
\def\ups{\upsilon}
\def\vp{\varphi}
\def\wh{\widehat}
\def\wt{\widetilde}
\def\x{\times}
\def\z{\zeta}
\def\({\left(}
\def\){\right)}
\def\[{\left[}
\def\]{\right]}
\def\<{\left<}
\def\>{\right>}

\def\SA{\mathcal A}
\def\SB{\mathcal B}
\def\SC{\mathcal C}
\def\SD{\mathcal D}
\def\SE{\mathcal E}
\def\SF{\mathcal F}
\def\SG{\mathcal G}
\def\SH{\mathcal H}
\def\SI{\mathcal I}
\def\SJ{\mathcal J}
\def\SK{\mathcal K}
\def\SL{\mathcal L}
\def\SM{\mathcal M}
\def\SN{\mathcal N}
\def\SO{\mathcal O}
\def\SP{\mathcal P}
\def\SQ{\mathcal Q}
\def\SR{\mathcal R}
\def\SS{\mathcal S}
\def\ST{\mathcal T}
\def\SU{\mathcal U}
\def\SV{\mathcal V}
\def\SW{\mathcal W}
\def\SX{\mathcal X}
\def\SY{\mathcal Y}
\def\SZ{\mathcal Z}

\newcommand{\BA}{\ensuremath{\mathbf A}}
\newcommand{\BB}{\ensuremath{\mathbf B}}
\newcommand{\BC}{\ensuremath{\mathbf C}}
\newcommand{\BD}{\ensuremath{\mathbf D}}
\newcommand{\BE}{\ensuremath{\mathbf E}}
\newcommand{\BF}{\ensuremath{\mathbf F}}
\newcommand{\BG}{\ensuremath{\mathbf G}}
\newcommand{\BH}{\ensuremath{\mathbf H}}
\newcommand{\BI}{\ensuremath{\mathbf I}}
\newcommand{\BJ}{\ensuremath{\mathbf J}}
\newcommand{\BK}{\ensuremath{\mathbf K}}
\newcommand{\BL}{\ensuremath{\mathbf L}}
\newcommand{\BM}{\ensuremath{\mathbf M}}
\newcommand{\BN}{\ensuremath{\mathbf N}}
\newcommand{\BO}{\ensuremath{\mathbf O}}
\newcommand{\BP}{\ensuremath{\mathbf P}}
\newcommand{\BQ}{\ensuremath{\mathbf Q}}
\newcommand{\BR}{\ensuremath{\mathbf R}}
\newcommand{\BS}{\ensuremath{\mathbf S}}
\newcommand{\BT}{\ensuremath{\mathbf T}}
\newcommand{\BU}{\ensuremath{\mathbf U}}
\newcommand{\BV}{\ensuremath{\mathbf V}}
\newcommand{\BW}{\ensuremath{\mathbf W}}
\newcommand{\BX}{\ensuremath{\mathbf X}}
\newcommand{\BY}{\ensuremath{\mathbf Y}}
\newcommand{\BZ}{\ensuremath{\mathbf Z}}

\def\bba{{\mathbb A}}
\def\bbb{{\mathbb B}}
\def\bbc{{\mathbb C}}
\def\bbd{{\mathbb D}}
\def\bbe{{\mathbb E}}
\def\bbf{{\mathbb F}}
\def\bbg{{\mathbb G}}
\def\bbh{{\mathbb H}}
\def\bbi{{\mathbb I}}
\def\bbj{{\mathbb J}}
\def\bbk{{\mathbb K}}
\def\bbl{{\mathbb L}}
\def\bbm{{\mathbb M}}
\def\bbn{{\mathbb N}}
\def\bbo{{\mathbb O}}
\def\bbp{{\mathbb P}}
\def\bbq{{\mathbb Q}}
\def\bbr{{\mathbb R}}
\def\bbs{{\mathbb S}}
\def\bbt{{\mathbb T}}
\def\bbu{{\mathbb U}}
\def\bbv{{\mathbb V}}
\def\bbw{{\mathbb W}}
\def\bbx{{\mathbb X}}
\def\bby{{\mathbb Y}}
\def\bbz{{\mathbb Z}}

\newcommand{\db}{\bar{\delta}}
\newcommand{\R}{\operatorname{r}}
\newcommand{\Z}{\mathbb{Z}}
\newcommand{\K}{\mathbb{K}}
\newcommand{\xp}{X^\prime}
\newcommand{\onto}{\twoheadrightarrow}
\newcommand{\Ext}{\operatorname{Ext}}
\newcommand{\red}{\textcolor{red}}

\newtheorem{thm}{Theorem}[section]
\newtheorem{lem}[thm]{Lemma}
\newtheorem{prop}[thm]{Proposition}
\newtheorem{cor}[thm]{Corollary}
\newtheorem*{1}{Stallings' Theorem (Integral)}
\newtheorem*{4.5}{Corollary \ref{free}}
\newtheorem*{4.1}{Theorem \ref{main}}
\newtheorem*{5.1}{Theorem \ref{universal}}
\newtheorem*{7}{Stallings' Theorem (Rational)}

\theoremstyle{definition}
\newtheorem{ex}[thm]{Example}
\newtheorem{defn}[thm]{Definition}

\newtheorem{rem}[thm]{Remark}

\newtheorem*{acknowledgement}{Acknowledgements}

\newcommand{\bz}{\mathbb{Z}}
\newcommand{\bq}{\mathbb{Q}}

\newcommand{\hz}[1]{\ensuremath{H_#1(-;\mathbb{Z})}}
\newcommand{\und}{\ensuremath{\un{\ \ }}}
\newcommand{\sll}{\ensuremath{S^3\backslash L}}
\newcommand{\sk}{\ensuremath{S^3\backslash K}}

\newcommand{\hq}[1]{\ensuremath{H_#1(-;\mathbb{Q})}}
\newcommand{\der}{\operatorname{Der}}
\newcommand{\id}{\operatorname{id}}
\newcommand{\ho}{\operatorname{Hom}}
\newcommand{\rank}{\operatorname{rank}}
\newcommand{\image}{\operatorname{image}}

\newcommand{\sss}{\scriptscriptstyle}
\newcommand{\2}[1]{\ensuremath{^{\sss (#1)}}}
\newcommand{\npp}{\ensuremath{^{\sss (n+1)}}}
\newcommand{\nm}{\ensuremath{^{\sss (n-1)}}}
\newcommand{\gn}{\ensuremath{G^{\sss (n)}_{\sss H}}}
\newcommand{\gnp}{\ensuremath{G^{\sss (n+1)}_{\sss H}}}
\newcommand{\gone}{\ensuremath{G^{\sss (1)}_{\sss H}}}
\newcommand{\an}{\ensuremath{A^{\sss (n)}_{\sss H}}}
\newcommand{\bn}{\ensuremath{B^{\sss (n)}_{\sss H}}}
\newcommand{\cn}{\ensuremath{C^{\sss (n)}_{\sss H}}}

\begin{document}

\title{Homology and derived series of groups}
\author{Tim Cochran\\Shelly Harvey}

\address{Department of Mathematics, Rice University\\Houston, TX 77005-1892, USA} 
\gtemail{\mailto{cochran@math.rice.edu}{\rm\qua 
and\qua}\mailto{shelly@math.rice.edu}}
\asciiemail{cochran@math.rice.edu, shelly@math.rice.edu} 

\begin{abstract} In 1964, John Stallings established an important
relationship between the low-dimensional homology of a group and
its lower central series. We establish a similar relationship
between the low-dimensional homology of a group and its
\emph{derived series}. We also define a \emph{torsion-free-solvable completion} of a group that is analogous
to the Malcev completion, with the role of the lower central series replaced by the derived series. We prove that the torsion-free-solvable completion is invariant under rational homology equivalence.
\end{abstract}

\asciiabstract{In 1964, John Stallings established an important
relationship between the low-dimensional homology of a group and its
lower central series.  We establish a similar relationship between the
low-dimensional homology of a group and its derived series. We also
define a torsion-free-solvable completion of a group that is analogous
to the Malcev completion, with the role of the lower central series
replaced by the derived series. We prove that the
torsion-free-solvable completion is invariant under rational homology
equivalence.}

\maketitle

\section{\label{intro}Introduction}

John Stallings, in his landmark paper \cite{St}, established the
following relationships between the low-dimensional homology of a
group and its lower central series. Recall that, for any ordinal
$\alpha$, the $\alpha^\supth$ term of the lower central series of
$G$, denoted $G_{\alpha}$, is inductively defined by $G_1=G$,
$G_{\alpha +1}=[G_{\alpha},G]$ and, if $\alpha$ is a limit
ordinal, $G_{\alpha}=\bigcap_{\beta<\alpha}G_{\beta}$. Stallings
also defines what we shall call the $\textsl{rational lower central
series}$, $G_{\alpha}^r$, which differs only in that $G_{\alpha
+1}^r$ consists of all those elements \emph{some finite power of
which} lies in $[G_{\alpha}^r,G]$. It is the most rapidly
descending central series whose successive quotients are torsion
free abelian groups.

\begin{1} {\rm\cite[Theorem 3.4]{St}}\qua Let $\phi\co A\to B$ be a
homomorphism that induces an isomorphism on $\hz1$ and an
epimorphism on $\hz2$. Then, for any finite $n$, $\phi$ induces
an isomorphism $A/A_n\cong B/B_n$. For the first infinite ordinal
$\omega$, it induces an embedding $A/A_{\omega}\subset
B/B_{\omega}$. If, in addition, $\phi$ is onto then, for each
ordinal $\alpha$, $\phi$ induces an isomorphism $A/A_{\alpha}\cong
B/B_{\alpha}$.
\end{1}

\begin{7} {\rm\cite[Theorem 7.3]{St}}\qua Let $\phi\co A\to B$ be a
homomorphism that induces an isomorphism on $\hq1$ and an
epimorphism on $\hq2$. Then, for all $\alpha\leq \omega$, $\phi$
induces an embedding $A/A^r_{\alpha}\subset B/B^r_{\alpha}$ and
for any finite $n$ it induces isomorphisms
$(A^r_n/A^r_{n+1})\otimes \mathbb{Q}\cong (B^r_n/B^r_{n+1})\otimes
\mathbb{Q}$.
\end{7}

These theorems have proven to be quite useful in topology. For
example if $A=\pi_1(\sll)$ where $L$ is a link of circles in
$S^3$, Stallings showed that the isomorphism type of each of the
quotients $A/A_n$ is an invariant of link concordance (even of
$I$--equivalence). The concordance invariance of Milnor's
$\overline{\mu}$--invariants was established by this means \cite{Ca}.
Stallings' theorems also give a criterion for establishing that a
collection of elements of a group generates a \emph{free}
subgroup.

Attempts have been made, most notably by Ralph Strebel \cite{Str},
to find a similar relationship between homology and the {\it
derived series} of groups, with limited success. We will use the
work of Strebel in a crucial way. Recall that the $\alpha^\supth$
term of the derived series, $G^{\sss (\alpha)}$, is defined by
$G\20=G$, $G^{\sss (\alpha +1)}=[G^{\sss (\alpha)},G^{\sss (\alpha)}]$ and, if
$\alpha$ is a limit ordinal,
$G^{\sss (\alpha)}=\bigcap_{\beta<\alpha}G^{\sss (\beta)}$. The derived
series has recently appeared prominently in joint work of the
first author, Kent Orr and Peter Teichner \cite{COT} \cite{COT2}
\cite{CT} \cite{T} where it was used to define new invariants for
classical knot concordance, and in work by other authors \cite{Ki1} \cite{Ki2}. It was
was also used to define \emph{higher-order Alexander invariants}
for knots \cite{C}, and 3--manifolds
\cite{Ha1}, and to define invariants for link concordance and rational homology cobordism of manifolds \cite{Ha2}. It has
also appeared recently in connection with questions about the
virtual first Betti number of 3--manifolds \cite{Ro} and in the study of complements of hyperplane arrangements \cite{PS}.

We show that there is a strong relationship between the
low-dimensional homology of a group and its derived series. For
example, we have the following strict analogue of Stallings'
Rational Theorem in the case $A$ is a free group.

\begin{4.5}\label{4.5} Let  $F$ be a free group and $B$ be a finitely related group (has a presentation with a finite number of relations). Let $\phi\co F\to B$ be a
homomorphism that induces
a monomorphism on $H_1(- ;\mathbb{Q})$ and an
epimorphism on $H_2(- ;\mathbb{Q})$. Then, for all $\alpha\leq \omega$, $\phi$
induces an embedding $F/F^{(\alpha)}\subset B/B^{(\alpha)}$.

\end{4.5}

The following examples show that several obvious generalizations
of Corollary~\ref{free} are false. Let $K$ be any knot in $S^3$,
$A=\pi_1(\sk)$ and $B=\bz$. Then the abelianization map yields a homomorphism
$\phi \co A\to A/[A,A]\cong\bz=B$. Then $\phi$ induces an isomorphism on all
integral homology groups (since, by Alexander Duality, $\sk$ has the
homology of a circle) and $B\2n=\{e\}$ for any $n\ge1$; whereas $A/A\2n$ is
known to be very large as long as the Alexander polynomial of $K$ is not $1$
\cite[Corollary 4.8]{C}. Thus $\phi$ cannot induce monomorphisms in general as in
Corollary~\ref{free}. Moreover, swapping the roles of $A$ and $B$ and
choosing a map $\bz\to A$ inducing an isomorphism on $\hz1$ gives a
homomorphism that again induces isomorphisms on all integral homology groups
and yet induces the map $\bz\to A/A\2n$ for each $n$, this being far
from surjective. Thus a direct analogue of Stallings' theorem might seem
hopeless.

However the second author, in a search for new invariants of link
concordance, introduced a new characteristic series,
$G\2n_{\sss H}$, associated to the derived series, called the
\emph{torsion-free derived series} (see \cite[Section 2]{Ha2}). Although this
series is not fully invariant, we show that it is functorial when
morphisms are restricted to those that induce a monomorphism on
$H_1(-;\mathbb{Q})$ and an epimorphism on $H_1(-;\mathbb{Q})$.
Using this we are able show the following analogue of Stallings'
Rational Theorem.

\begin{4.1}\label{4.1} Let $A$ be finitely-generated and $B$ finitely related (has a presentation with a finite number of relations). Let $\phi\co A\to B$ be a
homomorphism that induces a monomorphism on $H_1(- ;\mathbb{Q})$ and an
epimorphism on $H_2(- ;\mathbb{Q})$. Then, for each
$n\leq \omega$, $\phi$ induces a monomorphism
$A/A^{\sss (n)}_{\sss H} \subset B/B^{\sss (n)}_{\sss H}$. Moreover, if $\phi$
induces an isomorphism on $H_1(- ;\mathbb{Q})$ then, for each finite $n$,
$\an/A\npp_{\sss H} \to \bn/B\npp_{\sss H}$ is a monomorphism between modules of the
same rank (over $\bz[A/\an]$ and $\bz[B/\bn]$, respectively). If, in
addition $\phi$ is onto, then for each $n\leq \omega$, it induces an
isomorphism $A/\an \cong B/\bn$.
\end{4.1}

Corollary~\ref{free} is a special case of Theorem~\ref{main}.

Note that, by Stallings' Integral Theorem, the \emph{(pro-)nilpotent completion} of a group $G$,  $\varprojlim (G/G_n)$ is an invariant of integral homology equivalence. But neither the nilpotent completion nor its torsion-free version $\varprojlim (G/G^r_n)$ is invariant under \emph{rational} homology equivalence (for example the rational homology equivalence $\mathbb{Z}\to \mathbb{Z}$ where $t\to t^2$ does not induce an isomorphism on these completions). However there does exist a further completion, the \emph{Malcev completion} of $G$, which \emph{is} invariant under rational homology equivalence (a consequence of the last part of Stallings' Rational
Theorem). Turning to the derived series, the examples above show that neither the \emph{pro-solvable completion}, $\varprojlim (G/G^{\sss (n)})$, nor its torsion-free version, $\varprojlim (G/G^{\sss (n)}_r)$, nor even our $\varprojlim (G/G^{\sss (n)}_{\sss H})$ (see Section~\ref{basics}) is invariant under rational (or even integral) homology equivalence. However, early versions of \cite{Ha2} suggested another functor on groups, $G\to \wt G$, that we call the \emph{torsion-free-solvable completion} of $G$. This is properly viewed as an analogue of the Malcev completion of $G$, with the lower central series being replaced by the torsion-free derived series $G^{\sss (n)}_{\sss H}$. Here we rigorously construct $G\to \wt G$, whose kernel is $G^{\sss (\omega)}_{\sss H}$, and show that it is a rational homological localization functor, that is we show:

\begin{5.1}\label{5.1} Suppose $A$ is finitely generated and $B$ is finitely related. Suppose $\phi\co A\to
B$ induces an isomorphism on $H_1(-;\mathbb{Q})$ and an
epimorphism on $H_2(-;\mathbb{Q})$. Then $\phi$
induces an isomorphism $\wt A\cong \wt B$.
\end{5.1}

We also give the precise relationship between the torsion-free-solvable completion and the universal (integral) homological localization functor due to J Levine and P Vogel, and indicate the relation (by analogy) between the torsion-free-solvable completion and the Malcev completion.

\begin{acknowledgement}
The first author was partially supported by NSF DMS-0406573. The second author was partially supported by an NSF Postdoctoral Fellowship and NSF DMS-0539044.
\end{acknowledgement}

\section{The torsion-free derived series}\label{basics}

In this section we motivate and define Harvey's version of the
derived series, $G^{\sss (n)}_{\sss H}$, and remark on some of its elementary
properties.

If $G$ is a group then $G/G^{\sss (1)}$ is an abelian group but may have
$\bz$--torsion. If one would like to avoid this $\bz$--torsion then, in direct
analogy to the rational lower-central series above, one can define
$G\21_r=\{x\mid\exists k\neq0$ $x^k\in G\21\}$, slightly larger than $G\21$,
so that $G/G\21_r$ is $\bz$--torsion-free. Proceeding in this way,
defining $G\npp_r$ to be the radical of $[G\2n_r,G\2n_r]$, leads to
what has been called the \textsl{rational derived series of $G$}
\cite{Ha1, C, CT}. This is the most rapidly
descending series for which the quotients of successive terms are
$\bz$--torsion-free abelian groups. In \cite{Ha2} Harvey
observes that if a subgroup $G\2n_{\sss H}$ (normal in $G$) has been defined then
${G\2n_{\sss H}}/{[G\2n_{\sss H},G\2n_{\sss H}]}$ is not only an abelian group but also a right {\it module} over $\bz[G/G\2n_{\sss H}]$, where the action is induced from the
conjugation action of $G$ ($[x]g=[g^{-1}xg]$). One might seek to eliminate
torsion {\it in the module sense} from the successive quotients. This
motivated her definition of the $\textsl{torsion-free derived series}$ as
follows. Set $G\20_{\sss H}=G$. For $n\geq 0$, suppose inductively that
$G^{\sss (n)}_{\sss H}$ has been defined and is normal in $G$. Let $T_n$ be the subset
of ${\gn}/{[\gn,\gn]}$ consisting of the $\bz[G/G\2n_{\sss H}]$--torsion elements,
ie, the elements $[x]$ for which there exists some non-zero
$\g\in\bz[G/G\2n_{\sss H}]$, such that $[x]\g=0$. (In fact, since it will be
(inductively) shown below that $\bz [G/G\2n_{\sss H}]$ is an \emph{Ore Domain},
$T_n$ is a \emph{submodule}). Now consider the
epimorphism of groups:
$$
G^{\sss (n)}_{\sss H} \xrightarrow{\pi_n} \f{\gn}{[\gn,\gn]}
$$
and define $G^{\sss (n+1)}_{\sss H}$ to be the inverse image of $T_n$ under
$\pi_n$. Then $G^{\sss (n+1)}_{\sss H}$ is, by definition, a normal subgroup
of $G\2n_{\sss H}$ that contains $[\gn,\gn]$. It follows inductively that
$G^{\sss (n+1)}_{\sss H}$ contains $G^{\sss (n+1)}$ (and $G^{\sss (n+1)}_r$). Moreover,
since ${\gn}/{G\npp_{\sss H}}$ is the quotient of the module ${\gn}/{[\gn,\gn]}$
by its
torsion submodule, it is a $\bz[G/G^{\sss (n)}_{\sss H}]$ torsion-free module
\cite[Lemma 3.4]{Ste}.
Hence the successive quotients of the
torsion-free derived subgroups are {\it torsion-free modules} over
the appropriate rings. We define $G^{\sss (\omega)}_{\sss H}=\bigcap_{n<\omega}G_{\sss H}^{\sss (n)}$ as usual. However a discussion of
$G^{\sss (\alpha)}_{\sss H}$ for $\alpha >\omega$ would require further work since $\mathbb{Z}[G/G^{\sss (\omega)}_{\sss H}]$ is not in general an Ore Domain. We do not address this here.

Most of the following elementary properties
of the torsion-free derived series were established in \cite{Ha2}.
We repeat the proofs for the convenience of the reader.

\begin{prop}\label{PTFA} {\rm\cite{Ha2}}\qua $G/G^{\sss (n)}_{\sss H}$ is a \emph{poly-(torsion-free abelian) group}
(hereafter abbreviated PTFA), and consequently $\bz[G/G^{\sss (n)}_{\sss H}]$
is an Ore domain.
\end{prop}

Recall that a \textsl{poly-(torsion-free abelian) group} is one with
a finite subnormal series whose successive quotients are
torsion-free abelian groups. Such a group is solvable, torsion
free and locally indicable \cite[Proposition 1.9]{Str}.
Consequently $\bz[G/G^{\sss (n)}_{\sss H}]$ is an Ore domain and thus admits a
classical (right) division ring of fractions,
$\mathcal{K}(G/G^{\sss (n)}_{\sss H})$, into which it
embeds\cite[pages 591--592]{P}\cite [page 611]{P} \cite{Le}. Hence any
(right) module $M$ over $\bz[G/G^{\sss (n)}_{\sss H}]$ has a well-defined rank
which is defined to be the rank of the vector space $M
\otimes_{\bz[G/G^{\sss (n)}_{\sss H}]}\mathcal{K}(G/G^{\sss (n)}_{\sss H})$
\cite[page 48]{Co}. Alternatively the rank can be defined to be the
maximal integer $m$ such that $M$ contains a submodule isomorphic
to $(\bz[G/G^{\sss (n)}_{\sss H}])^m$.

\begin{prop}\label{normality} {\rm\cite{Ha2}}\qua $G^{\sss (n+1)}_{\sss H}$ is a normal subgroup
of $G$.
\end{prop}

\begin{proof}[Proof of Proposition~\ref{normality}] The proof is by
induction on $n$. We assume $G^{\sss (n)}_{\sss H}$ is normal in $G$ and show that $G^{\sss (n+1)}_{\sss H}$
is normal in $G$. Consider $x\in G^{\sss (n+1)}_{\sss H}$ and $g\in G$. Then $x$ and
$g^{-1}xg$ lie in $G^{\sss (n)}_{\sss H}$ and, by definition of the module structure on
${\gn}/{[\gn,\gn]}$, $\pi _n(g^{-1}xg)=\pi _n(x)g$. Since $x\in
G^{\sss (n+1)}_{\sss H}$,  $\pi _n(x)$ is torsion. To show that $g^{-1}xg\in
G^{\sss (n+1)}_{\sss H}$, we need to show that $\pi _n(x)g$ is torsion. Hence the
normality of $G^{\sss (n+1)}_{\sss H}$ will follow from showing that the torsion subgroup
$T_n$ is a \emph{submodule}. But the set of torsion elements of any module over an
Ore domain is known to be a submodule \cite[page 57]{Ste}. Thus the desired
result follows from Proposition~\ref{PTFA}.
\end{proof}

The torsion-free derived subgroups are characteristic subgroups,
but are not \emph{fully invariant}, that is an arbitrary
homomorphism of groups $\phi \co A \lra B$ need not send $A^{\sss (n)}_{\sss H}$
into $B^{\sss (n)}_{\sss H}$. The simplest example is $A=\langle x,y,z| [z,[x,y]]=1 \rangle$, $B=\langle x,y\rangle$, with $\phi$ sending $(x\to x,y\to y, z\to 0)$. Then $[x,y]\in A^{\sss (2)}_{\sss H}$ but $\phi ([x,y])$ is not in $B^{\sss (2)}_{\sss H}=B^{\sss (2)}$. However, the following observation is used in
our main theorem to show (inductively) that if $\phi \co A\to B$ induces a
monomorphism on $H_1(- ;\mathbb{Q})$ and an epimorphism on $H_2(-
;\mathbb{Q})$, and if $A$ is finitely generated and $B$ finitely
related, then $\phi$ \emph{does} send $A^{\sss (n)}_{\sss H}$ into
$B^{\sss (n)}_{\sss H}$. It follows that the torsion-free derived series is
functorial in this restricted category.

\begin{prop}\label{functoriality} {\rm\cite{Ha2}}\qua If $\phi \co A \to B$ induces a
monomorphism $\phi \co A/A^{\sss (n)}_{\sss H} \to B/B^{\sss (n)}_{\sss H}$, then $\phi
(A^{\sss (n+1)}_{\sss H})\subset B^{\sss (n+1)}_{\sss H}$ and hence $\phi$ induces a homomorphism
$\phi :$\break $A/A^{\sss (n+1)}_{\sss H} \to B/B^{\sss (n+1)}_{\sss H}$.
\end{prop}

\begin{proof}[Proof of Proposition~\ref{functoriality}] Note that the
hypothesis implies that $\phi$ induces a ring monomorphism
$\tilde{\phi}\co \bz[A/A^{\sss (n)}_{\sss H}]\to\bz[B/B^{\sss (n)}_{\sss H}]$. Suppose
that $x\in A^{\sss (n+1)}_{\sss H}$. Consider the diagram below.
\begin{equation*}
\begin{CD}
A^{\sss (n)}_{\sss H}          @>\pi_A>>   \f{\an}{A\npp_{\sss H}}\\
@VV\phi V             @VV\bar{\phi} V\\
B^{\sss (n)}_{\sss H}          @>\pi_B>>   \f{\bn}{B\npp_{\sss H}}\\
\end{CD}
\end{equation*}
By definition, $\pi_A(x)$ is torsion, that is there is some non-zero
$\gamma\in\bz[A/A^{\sss (n)}_{\sss H}]$ such that $(\pi_A(x))\gamma =0$. It is easy to
check that $\bar{\phi}$ is a homomorphism of right
$\bz[A/A^{\sss (n)}_{\sss H}]$--modules using the module structure induced on
${\bn}/{B\npp_{\sss H}}$ by $\tilde{\phi}$ (since $\phi (a^{-1}xa)=\phi(a)^{-1}
\phi(x) \phi(a)$). Thus $(\bar{\phi}(\pi_A(x)))\tilde\phi (\gamma)=0$ and so, since $\tilde{\phi}$ is injective,
$\bar{\phi}(\pi_A(x))$ is a $\bz[B/B^{\sss (n)}_{\sss H}]$--torsion element. But
$\bar{\phi}(\pi_A(x))=\pi_B(\phi (x))$, showing that $\phi(x)\in
B^{\sss (n+1)}_{\sss H}$. Hence $\phi(A^{\sss (n+1)}_{\sss H})\subset B^{\sss (n+1)}_{\sss H}$.
 \end{proof}

For some groups, such as free groups and free-solvable groups $F/F^{\sss (k)}$, the derived series and the torsion-free derived series
coincide.

\begin{prop}\label{torsionfree} {\rm\cite{Ha2}}\qua If $G$ is a group such that, for every $n$,
$G^{\sss (n)}/G^{\sss (n+1)}$ is a $\bz[G/G^{\sss (n)}]$--torsion-free
module, then the torsion-free derived series of $G$ agrees with
the derived series of $G$. Hence for a free group $F$, $F^{\sss (n)}_{\sss H}=F^{\sss (n)}$ for each $n$, and similarly for a free solvable group $F/F^{(k)}$.
\end{prop}

\begin{proof}[Proof of Proposition~\ref{torsionfree}] By definition, $G^{\sss (0)}_{\sss H}=G^{\sss (0)}=G$. Suppose
$G^{\sss (n)}_{\sss H}=G^{\sss (n)}$. Then, under
the hypotheses, $G^{\sss (n)}_{\sss H}/[G^{\sss (n)}_{\sss H},G^{\sss (n)}_{\sss H}]$ is a torsion-free module and so
$G^{\sss (n+1)}_{\sss H}=\ker\pi_n=[G^{\sss (n)}_{\sss H},G^{\sss (n)}_{\sss H}]=[G^{\sss (n)},G^{\sss (n)}]=G^{\sss (n+1)}$.

It is well known that $F^{\sss (n)}/F^{\sss (n+1)}$ is a $\bz [F/F^{\sss (n)}]$--torsion-free
module. This can be seen by examining the free $\bz[F/F^{\sss (n)}]$
cellular chain complex for the covering space of a wedge of circles
corresponding to the subgroup
 $F^{\sss (n)}$. The module $F^{\sss (n)}/F^{\sss (n+1)}$ is merely the first homology of
this chain complex. Since the chain complex can be chosen to have
no $2$--cells, its first homology
 is a submodule of a free module and thus is torsion-free. Hence the derived
series and the
 torsion-free derived series of a free group agree. The same is true for any free-solvable group $F/F^{(k)}$.
 \end{proof}
The following basic result is useful.

\begin{prop}\label{idempotency} If $\phi \co A \to B$ is surjective and $\ker\phi \subset A^{\sss (m)}_{\sss H}$ then $\phi$ induces
isomorphisms $A/A^{\sss (n)}_{\sss H} \cong B/B^{\sss (n)}_{\sss H}$ for all $n\leq m$. In particular, $(G/\gn)^{\sss (n)}_{\sss H}=0$.
\end{prop}
\begin{proof}[Proof of Proposition~\ref{idempotency}] By induction we assume $\phi$ induces isomorphisms $A/A^{\sss (i)}_{\sss H} \cong B/B^{\sss (i)}_{\sss H}$ for each $i\leq n$ for some $0\leq n<m$. Since $\phi$ induces a map $A/A^{\sss (n)}_{\sss H} \to B/B^{\sss (n)}_{\sss H}$ surely $\phi(A^{\sss (n)}_{\sss H})\subset B^{\sss (n)}_{\sss H}$. We claim that $\phi(A^{\sss (n)}_{\sss H})= B^{\sss (n)}_{\sss H}$, because for any $b\in B^{\sss (n)}_{\sss H}$ there is an $a\in A$ such that $\phi(a)=b$ and, by the inductive hypothesis, we must have $a\in A^{\sss (n)}_{\sss H}$. Therefore it also follows that $\phi([A^{\sss (n)}_{\sss H},A^{\sss (n)}_{\sss H}])= [B^{\sss (n)}_{\sss H},B^{\sss (n)}_{\sss H}]$, a fact we use below. Moreover, by Proposition~\ref{functoriality}, $\phi(A^{\sss (n+1)}_{\sss H})\subset B^{\sss (n+1)}_{\sss H}$ and so induces an epimorphism $\phi \co A/A^{\sss (n+1)}_{\sss H} \to B/B^{\sss (n+1)}_{\sss H}$. We need to show this is injective to complete the proof. Suppose there exists $a\in A$ such that $\phi (a)=b\in B^{\sss (n+1)}_{\sss H}$. By the inductive hypothesis, $a\in A^{\sss (n)}_{\sss H}$. Refer to the commutative diagram below.
\begin{equation*}
\begin{CD}
A^{\sss (n)}_{\sss H}          @>\pi_A>>   \f{\an}{[A^{\sss (n)}_{\sss H},A^{\sss (n)}_{\sss H}]}\\
@VV\phi V             @VV\bar{\phi} V\\
B^{\sss (n)}_{\sss H}          @>\pi_B>>   \f{\bn}{[B^{\sss (n)}_{\sss H},B^{\sss (n)}_{\sss H}]}\\
\end{CD}
\end{equation*}
Since  $b\in B^{\sss (n+1)}_{\sss H}$, $\pi_B(b)$ is a torsion element and so there exists a non-zero $\gamma \in \mathbb{Z}[B/\bn]$ such that $(\pi_B(b))\gamma=0$. We seek to express this entirely in terms of $B$ itself. Lift $\gamma$ to $\mathbb{Z}B$ and express it as $\sum k_i\eta_i$ where $k_i\in \mathbb{Z}$ and $\eta_i\in B$. The condition that $(\pi_B(b))\gamma=0$ translates to the condition $\prod (\eta^{-1}_ib^{k_i}\eta_i)=z$  for some $z\in [B^{\sss (n)}_{\sss H},B^{\sss (n)}_{\sss H}]$. Since $\phi$ is surjective, $\eta_i=\phi(\alpha_i)$ for some $\alpha_i\in A$ and $z=\phi(x)$ for some $x\in [A^{\sss (n)}_{\sss H},A^{\sss (n)}_{\sss H}]$ by our remark above. Then $x^{-1}\prod (\alpha^{-1}_ia^{k_i}\alpha_i)$ is in the kernel of $\phi$ and hence by hypothesis lies in $A^{\sss (m)}_{\sss H}\subset A^{\sss (n+1)}_{\sss H}$. Since $x\in [A^{\sss (n)}_{\sss H},A^{\sss (n)}_{\sss H}]\subset A^{\sss (n+1)}_{\sss H}$, it follows that $\prod (\alpha^{-1}_ia^{k_i}\alpha_i)$ is in $A^{\sss (n+1)}_{\sss H}$. Let $\beta=\sum k_i\alpha_i\in \mathbb{Z}[A/\an]$ which is non-zero since $\phi(\beta)=\gamma$ is non-zero. Note that $\pi_A(a)\beta$ is represented by $\prod (\alpha^{-1}_ia^{k_i}\alpha_i)$. Thus we have shown that $\pi_A(a)\beta$ is zero when mapped into $A^{\sss (n)}_{\sss H}/A^{\sss (n+1)}_{\sss H}$ and hence that the image of $\pi_A(a)$ is torsion in this module. But the successive quotients of terms in the torsion-free derived series are torsion-free modules as remarked earlier, so $\pi_A(a)=0$ in this module. Therefore $a\in A^{\sss (n+1)}_{\sss H}$, completing the proof of injectivity.

To prove that $(G/\gn)^{\sss (n)}_{\sss H}=0$, apply the above to the $\phi\co G\to G/\gn$.
\end{proof}
We give some elementary examples of the torsion-free derived series.

\begin{ex}\label{finite} Elements of finite order in a group $G$ are contained in every term of the rational derived series and hence every term of the torsion-free derived series. Therefore the torsion-free derived series of a finite group stabilizes at $n=0$, that is $\gn = G^{\sss (0)}_{\sss H}=G$. In general $G/G^{\sss (\omega)}_{\sss H}$ is a torsion-free group that is trivial if $G$ is finite. More generally, if $\beta _1(G)=0$, that is if $H_1(G;\mathbb{Q})=0$, then the torsion-free derived series, like the rational derived series, stabilizes at $n=0$.
\end{ex}
\begin{ex}\label{nilpotent} If $G$ is nilpotent with $\beta _1(G)>0$ then the terms of the torsion-free derived series are much larger than the terms of the rational derived series. In this case the torsion-free derived series stabilizes at $n=1$. For suppose $x\in G$ is an arbitrary element of $\gn$ for some $n\geq 1$ and $t\in G$ is non-zero in $G/G^{\sss (1)}_{\sss H}$. We will show that $x\in G^{\sss (n+1)}_{\sss H}$. Since G is nilpotent, for some $k$ the simple $(k+1)$--commutator, $[t,[t,[...[t,x]...]]]$ is the identity in $G$. In the module $\gn /[\gn ,\gn]$ this gives the relation $(t-1)^kx=0$, showing that $x$ represents a torsion element (since $(t-1)^k$ is non-zero in the integral domain $\mathbb{Z}[G/\gn ]$). Hence $x\in G^{\sss (n+1)}_{\sss H}$.
\end{ex}
\begin{ex}\label{freesolvable} If $G$ is free-solvable of class $n$, that is $G=F/F^{\sss (n)}$ where $F$ is a non-abelian free group, then by Proposition~\ref{torsionfree} and Proposition~\ref{idempotency}, $G/G^{\sss (m)}_{\sss H} \cong F/F^{\sss (m)}$ for all $m\leq n$ and $G^{\sss (n)}=0$. Thus the torsion-free derived series stabilizes at precisely $n$.
\end{ex}
\begin{ex}\label{bettione} If $\beta _1(G)=1$, for example if $G$ is
  the fundamental group of the exterior of a knotted circle in $S^3$,
  then the torsion-free derived series stabilizes at $n=1$ because the
  classical Alexander module, $G^{\sss (1)}_{\sss H}/[G^{\sss
      (1)}_{\sss H},G^{\sss (1)}_{\sss H}]$ can be seen to be a
  torsion module over $\mathbb{Z}[t,t^{-1}]$ (see, for example,
  \cite[Proposition 2.1]{COT}). Thus $G^{\sss (2)}_{\sss H}=G^{\sss (1)}_{\sss H}$. Alternatively, note that the epimorphism $G\to \mathbb{Z}$ induces an isomorphism on $H_1(-;\mathbb{Q})$ and an epimorphism on $H_2(-;\mathbb{Q})$ and so, by Theorem~\ref{main}, it induces isomorphisms $G/\gn \cong \mathbb{Z}$.
\end{ex}
\begin{ex}\label{rankzero} Generalizing Example~\ref{bettione}, if $G$ is any group whose classical Alexander module $G^{\sss (1)}_{\sss H}/[G^{\sss (1)}_{\sss H},G^{\sss (1)}_{\sss H}]$ has rank zero (over the multivariable Laurent polynomial ring $\mathbb{Z}[t_1^{\pm 1},...,t_m^{\pm1}]$) then the torsion-free derived series stabilizes at $n\leq 1$. The fundamental group of any $3$--manifold that fibers over a circle has this property \cite[Proposition 8.4]{Ha1}. By contrast, the derived series of the fundamental group of a knot exterior does not stabilize for any finite $n$ unless the Alexander polynomial is $1$ \cite[Corollary 4.8]{C}. More generally, if the rank of one of the \emph{higher-order modules} $\gn /[\gn ,\gn ]$, is zero then the torsion-free derived series stabilizes at the minimum value of $n$ for which this holds.
\end{ex}
\begin{ex}\label{general} For groups with $\beta _1(G)\geq 2$, a wide variety of behavior is possible. For example, as previously mentioned, the torsion-free derived series of a non-abelian free group, $F$, agrees with the derived series and is known to be highly non-trivial, stabilizing at $\omega$ (in fact $F^{\sss (\omega)}_{\sss H}=0$). There are many non-free groups where the torsion-free derived series does not stabilize at a finite ordinal. Recall that a link $\{ L_i\}$ of $m$ circles in $S^3$ is called a \emph{boundary link} if there is a collection of $m$ compact oriented surfaces $\{V_i\}$, disjointly embedded in $S^3$ such that $\partial V_i=L_i$. If $G=\pi_1(S^3-L)$ then the meridional map $F\to G$ has a right inverse $G\to F$ (use the Pontryagin construction to get a map from $S^3-L$ to a wedge of circles). The latter epimorphism induces an isomorphism on $H_1(-;\mathbb{Z})$ and an epimorphism $H_2(-;\mathbb{Z})$ and so, by Theorem~\ref{main}, for $n\leq \omega$, $G/\gn \cong F/F^{\sss (n)}$ and $G/G^{\sss (\omega)}_{\sss H}\cong F$. It follows from Stallings' Integral Theorem that $G^{\sss (\omega)}_{\sss H}=G_\omega$. Since for a boundary link it is known that $G_\omega/[G_\omega,G_\omega]$ is a torsion module, the torsion-free derived series stabilizes at $\omega$. However, suppose $G$ is the fundamental group of a $2$--complex $X$ with $H_1(X;\mathbb{Z})\cong \mathbb{Z}^m$ and $H_2(X;\mathbb{Z})=0$. Then there is a map $F\to G$ inducing isomorphisms on $H_1(-;\mathbb{Z})$ and $H_2(-;\mathbb{Z})$. Thus, by Theorem~\ref{main}, $F/F^{\sss (n)}\subset G/\gn$. However generally this is not surjective and $\gn / G^{\sss (n+1)}_{\sss H}$ is generally not even abstractly isomorphic to $F^{\sss (n)}/F^{\sss (n+1)}$ (see Remark~\ref{noniso}). In these cases it is not in general known if the torsion-free derived series stabilizes at $\omega$.
\end{ex}

It will be important for our main theorem that the reader understand the connection between the torsion-free derived series and group homology. This is provided by the following basic observations.

\begin{rem}\label{observation} $\an/[\an,\an]\cong H_1(A;\mathbb{Z}[A/\an])$. For an algebraist this is a consequence of the definition $H_1(A;\mathbb{Z}[A/\an])\equiv Tor_1^{A}(\mathbb{Z}[A/\an],\mathbb{Z})$ and the easy observation that the latter is $Tor_1^{\an}(\mathbb{Z},\mathbb{Z})\cong \an/[\an,\an]$ \cite[Lemma 6.2]{HS}. For a topologist, $H_1(A;\mathbb{Z}[A/\an]
)$ is thought of as the first homology with twisted coefficients of an
  aspherical space $K(A,1)$ where $\pi_1(K(A,1))\cong A$ and the
  coefficient system is induced by $\pi_1(K(A,1))\cong A\to A/\an$
  \cite[page 335]{HS}. Then $H_1(K(A,1);\mathbb{Z}[A/\an])$ can be interpreted as the first homology module of the covering space of K(A,1) corresponding to the subgroup $\an$, which is $\an/[\an,\an]$ \cite[Theorems VI3.4 and 3.4*]{Wh}.
 \end{rem}

\begin{prop}\label{observation2}$\phantom{99}$
\begin{enumerate} 
\item $\an/A^{\sss (n+1)}_{\sss H}$ is equal to $H_1(A;\mathbb{Z}[A/\an])$ modulo its $\mathbb{Z}[A/\an]$--torsion submodule.
\item  $A^{\sss (n+1)}_{\sss H}$ is the kernel of the composition:
\begin{align*}A^{\sss (n)}_{\sss H} \xrightarrow{\pi_n}
  \f{\an}{[\an,\an]}&=H_1(A;\mathbb{Z}[{A}/{\an}])\\
&\to H_1(A;\mathbb{Z}[A/\an])\otimes_{\mathbb{Z}[A/\an]}\SK (A/\an).
\end{align*}
\end{enumerate}\end{prop}

\begin{proof} Property $1)$ follows directly from
  Remark~\ref{observation} and the definition of $A^{\sss (n+1)}_{\sss
    H}$. For Property $2)$, note that tensoring with the quotient
  field $\SK (A/\an)$ kills precisely the $\mathbb{Z}[A/\an]$--torsion
  submodule \cite[Corollary II.3.3]{Ste}.
\end{proof}

 \section{Topological applications}

Our theorems, like those of Stallings, have applications to the
study of rational homology equivalences of topological spaces and to the
study of links in particular. We also give an application to when
 a set of elements of a group generates a free subgroup (similar to that of
Stallings). Deeper secondary applications will appear in \cite{Ha2}.

\begin{thm}\label{topological}Let $f\co X\to Y$ be a continuous
map between connected CW complexes that induces a monomorphism on
$H_1(- ;\mathbb{Q})$ and an epimorphism on $H_2(- ;\mathbb{Q})$.
Let $A$, $B$  denote  $\pi_1(X)$, $\pi_1(Y)$ respectively. Suppose
that $A$ is finitely-generated and $B$ is finitely related. Then,
for each $n\leq \omega$, $f$ induces a monomorphism $A/\an \subset
B/\bn$. If, additionally, $f$ induces an isomorphism on $H_1(-
;\mathbb{Q})$ then, for any finite $n$,
$\an/A\npp_{\sss H}\to\bn/B\npp_{\sss H}$ is a monomorphism between modules
of the same rank (over $\mathbb{Z}[A/\an]$ and
$\mathbb{Z}[B/\bn]$, respectively). If, in addition
$f_*\co \pi_1(X)\to \pi_1(Y)$ is onto, then, for each finite $n$, $f$
induces an isomorphism $A/\an \cong B/\bn$.
\end{thm}

\begin{proof} The theorem follows immediately from
Theorem~\ref{main} once we make the well-known observation that
 since $f$ induces an epimorphism $H_2(X;\mathbb{Q})\rightarrow H_2(Y;\mathbb{Q})$, it induces
 an epimorphism $H_2(A;\mathbb{Q})\rightarrow
 H_2(B;\mathbb{Q})$.
\end{proof}
\begin{cor}\label{homologycobordism} If $X$ and $Y$ are $n$--manifolds
(possibly with boundary) that are rationally homology cobordant
(relative their boundary) via the $(n+1)$--manifold $Z$, and $A$,
$B$ and $C$ denote their respective fundamental groups, then for
any n, the inclusion maps induce monomorphisms $A/\an \subset
C/\cn$, and $B/\bn \subset C/\cn$. Moreover $\an/A\npp_{\sss H}$ and
$\bn/B\npp_{\sss H}$ are modules of the same rank (over
$\mathbb{Z}[A/\an]$ and $\mathbb{Z}[B/\bn]$, respectively).
\end{cor}

In \cite{Ha2}, Harvey uses
Corollary~\ref{homologycobordism} to give other new
invariants of homology cobordism of $(2k-1)$--manifolds, using the
von Neumann $\rho$ --invariants of J Cheeger and M Gromov. In addition, the ranks of the above modules are new homology cobordism invariants. A
common example of this Corollary is the case of the exteriors of
concordant links.

\begin{cor}\label{Iequivalence} If $L_0$ and $L_1$ are
 compact subsets of $S^n$ that are concordant (or even merely $I$--equivalent) via the subset $W$
 of $S^n\times [0,1]$ and if $A$,$B$ and $C$ denote
the fundamental groups of the exteriors $S^n-L_0$, $S^n-L_1$ and
$S^n\times [0,1]-W$ respectively then, for any n, the inclusion
maps induce monomorphisms $A/\an \subset C/\cn$, and $B/\bn
\subset C/\cn$. Moreover $\an/A\npp_{\sss H}$ and $\bn/B\npp_{\sss H}$ are
modules of the same rank (over $\mathbb{Z}[A/\an]$ and
$\mathbb{Z}[B/\bn]$, respectively).
\end{cor}

\begin{proof} By Alexander Duality the exterior $S^n\times [0,1]-W$
is a homology cobordism (or, in the case of $I$--equivalence, just a
homology product) between the exteriors of $L_0$ and $L_1$ and so
Theorem~\ref{topological} applies.
\end{proof}

As above, Harvey uses Corollary~\ref{Iequivalence} to
give other new results in link concordance \cite{Ha2}. For example, note that the ranks of the above modules are new concordance invariants of links, generalizing the well-known fact that the
rank of the Alexander module of a link is an invariant of
concordance. It was this particular application that
first motivated Harvey's definition of the torsion-free derived
series. The above-mentioned Cheeger--Gromov invariants are used in \cite{Ha2} to show that the concordance group of disk links in any odd dimension has infinite rank even modulo local knotting.

Stallings theorem also gives a beautiful criterion to
establish that a set of elements of a group generates a free
subgroup. We have our own version, although we do not know an
example where it is stronger than Stallings result and there are
examples where it is weaker.

\begin{prop}\label{generatesfree} If $G$ is a finitely-presented
group with $H_2(G;\mathbb{Q})=0$, and ${x_1,...,x_k}$ is a set of
elements of $G$ that is linearly independent in the abelianization
of $G$ then this set generates a free subgroup in
$G/G^{\sss (\omega)}_{\sss H}$.
\end{prop}
\begin{proof} Following Stallings, let $F$ be a free group of rank
$k$ equipped with the obvious map into $G$ determined by the
${x_i}$. Then apply Corollary~\ref{free}. The result then follows
easily since for the free group $F^{\sss (\omega)}_{\sss H}\cong
F^{\sss (\omega)}=0$.
\end{proof}

\section{Proof of the Main Theorem}

In this section we prove our main result, Theorem~\ref{main},
which is a direct analogue, for the torsion-free derived series,
of Stallings' Rational Theorem for the rational lower
central series.

\begin{thm}\label{main}  Let $A$ be finitely-generated and $B$ finitely related. Suppose $\phi\co A\to B$ induces
a monomorphism on $H_1(- ;\mathbb{Q})$ and an
epimorphism on $H_2(- ;\mathbb{Q})$. Then, for each
$n\leq \omega$, $\phi$ induces a monomorphism $A/\an \subset
B/\bn$. Moreover, if $\phi$ induces an isomorphism on
$H_1(- ;\mathbb{Q})$ then, for any finite $n$, then
$\an/A\npp_{\sss H}\to\bn/B\npp_{\sss H}$ is a monomorphism between modules
of the same rank (over $\mathbb{Z}[A/\an]$ and
$\mathbb{Z}[B/\bn]$, respectively). If, in addition, $\phi$ is
onto, then for each $n\leq \omega$ it induces an isomorphism $A/\an
\cong B/\bn$. In the case that $(B,A)$ admits the structure of a relative 2--complex, the first and third conclusions above remain valid without the finiteness assumptions on $A$ and $B$.
\end{thm}

\begin{proof}[Proof of Theorem \ref{main}] The proof of the first claim is by induction on
$n$. The case $n=1$ is clear since $A/A^{\sss (1)}_{\sss H}$ is merely
$H_1(A;\mathbb{Z})/\{\mathbb{Z}\text{--Torsion}\}$ and the hypothesis that $\phi$ induces
a monomorphism on $H_1(-;\mathbb{Q})$ implies that it also
induces a monomorphism on $H_1(-;\mathbb{Z})$ modulo torsion.
Now assume that the first claim holds for $n$, ie, $\phi$ induces a monomorphism $A/\an \subset
B/\bn$. We will prove that it
holds for $n+1$.

It follows from Proposition \ref{functoriality} that $\phi(A\npp_{\sss H})\subset B\npp_{\sss H}$. Hence the diagram below exists and is commutative. In light of the Five Lemma, we see that it suffices to show that
$\phi$ induces a monomorphism $\an/A\npp_{\sss H}\to\bn/B\npp_{\sss H}$.
\begin{equation*}
\begin{CD}
1      @>>>    \f{\an}{A\npp_{\sss H}}  @>>>   \f A{A\npp_{\sss H}}  @>>>
\f A{\an}    @>>> 1\\
&&     @VV \phi V        @VV\phi_{n+1}V       @VV\phi_nV\\
1      @>>>    \f{\bn}{B\npp_{\sss H}}  @>>>   \f B{B\npp_{\sss H}}  @>>>
\f B{\bn}    @>>> 1\\
\end{CD}
\end{equation*}
For simplicity we abbreviate $A/\an$ by $A_n$ and $B/\bn$ by
$B_n$. The inductive hypothesis is that $\phi$ induces a
monomorphism $A_n\to B_n$ and hence a ring monomorphism
$\mathbb{Z} A_n\to\mathbb{Z} B_n$. Since $A_n$ and $B_n$ are PTFA
groups by Proposition~\ref{PTFA}, the rings $\mathbb{Z} A_n$ and
$\mathbb{Z} B_n$ are right Ore domains and so admit classical right
rings of quotients $\SK(A_n)$ and $\SK(B_n)$, respectively. Hence
$\phi$ induces a monomorphism $\SK(A_n)\to\SK(B_n)$, which endows
$\SK(B_n)$ with the structure of a $\SK(A_n)-\SK(B_n)$ bimodule.

Suppose that $\an/A\npp_{\sss H}\to\bn/B\npp_{\sss H}$ were \emph{not} injective. By examining the diagram below, we see that there exists an $a\in\an$ representing a \emph{non}-torsion class $[a]$ in ${\an}/{[\an,\an]}$ such that $\phi (a)$
represents a $\mathbb{Z} B_n$--torsion class, $[\phi(a)]$, in ${\bn}/{[\bn,\bn]}$.
\begin{equation*}
\begin{CD}
\f{\an}{[\an,\an]}          @>\pi_A>>   \f{\an}{A\npp_{\sss H}}\\
@VV\phi V             @VV\phi V\\
\f{\bn}{[\bn,\bn]}          @>\pi_B>>   \f{\bn}{B\npp_{\sss H}}\\
\end{CD}
\end{equation*}
Now consider the following commutative diagram where the horizontal equivalences follow from Remark~\ref{observation}.
\begin{equation*}
\begin{CD}
\f{\an}{[\an,\an]}\\
@VV i V \\
\f{\an}{[\an,\an]}\ox_{\mathbb{Z} A_n}\SK A_n\\
@VV\id\ox\phi V\\
\f{\an}{[\an,\an]}\ox_{\mathbb{Z} A_n}\SK B_n @>\cong>>
H_1(A;\mathbb{Z} A_n)\ox_{\mathbb{Z} A_n}\SK B_n\\
@VV\phi\ox\id V  @VV\phi\ox\id V\\
\f{\bn}{[\bn,\bn]}\ox_{\mathbb{Z} B_n}\SK B_n @>\cong>>
H_1(B;\mathbb{Z} B_n)\ox_{\mathbb{Z} B_n}\SK B_n
\end{CD}
\end{equation*}
We assert (and shall establish below) that the kernel of the vertical composition
 $\psi=(\phi\ox\id)\circ(\id\ox\phi)\circ i$ is the $\mathbb{Z} A_n$--torsion submodule of $\an/[\an ,\an ]$.
Assuming this, we finish the inductive proof.
Note that, $\psi([a])$ is $[\phi(a)]\ox1$. Since $[\phi(a)]$ is
 $\mathbb{Z} B_n$--torsion, $\psi([a])=0$, since tensoring with the
 quotient field kills the torsion submodule \cite[Corollary II.3.3]{Ste}. By the assertion, $[a]$ is $\mathbb{Z} A_n$--torsion, contradicting the choice of $[a]$. This contradiction shows that $\an/A\npp_{\sss H}\to\bn/B\npp_{\sss H}$ is
injective. The injectivity at the first infinite ordinal follows immediately, finishing the inductive step of the proof of the first
part of Theorem~\ref{main}, modulo our assertion.

Now we set out to establish that the kernel of $\psi$ is the $\mathbb{Z} A_n$--torsion submodule of $\an/[\an ,\an ]$. The kernel of $i$ is precisely this submodule since tensoring with the quotient field kills precisely the torsion submodule. Therefore it suffices to show that the other two maps in the composition are injective.

First note that $\id\ox\phi$ is injective by an
application of the following Lemma~\ref{mono}.

\begin{lem}\label{mono} Suppose $H\subset G$ and $\mathbb{Z} H$ and $\mathbb{Z} G$
are Ore Domains. Then for any right $\mathbb{Z} H$--module $M$, the
map
$$
id\ox i\co  M\ox_{\mathbb{Z} H}\SK H\to M\ox_{\mathbb{Z} H}\SK G
$$
 is a monomorphism (of right $\SK H$--modules).
Moreover, the $\mathcal{K}H$--rank of the domain of this map equals
the $\mathcal{K}G$--rank of the range.
\end{lem}

\begin{proof}[Proof of Lemma~\ref{mono}] Since $\SK H$ is a $\mathbb{Z}
H-\SK H$ bimodule, $M\ox_{\mathbb{Z} H}\SK H$ is a right $\SK
H$--module. Since any $\SK H$--module is free
\cite[Proposition I.2.3]{Ste}, there is some index set $I$ such that
$\psi\co M\ox_{\mathbb{Z} H}\SK H\cong\op_I\SK H$ as right $\SK
H$--modules. Thus the $\mathcal{K}H$--rank of the domain of $id\ox
i$ is the cardinality of $I$. Hence
$$
\psi\ox\id\co  \(M\otimes_{\mathbb{Z} H}\SK H\)\otimes_{\SK H}\SK
G\lra\(\oplus_I\SK H\)\otimes_{\SK H}\SK G
$$
is an isomorphism of right $\SK G$--modules, and hence of right
$\SK H$--modules. Since the domain of $\psi\ox\id$ is isomorphic to
$M\ox_{\mathbb{Z} H}\SK G$, we can use $\psi$ and $\psi\ox\id$ to
see that the first part of the Lemma is equivalent to showing that
$\op_I\SK H\lra(\op_I\SK H)\ox_{\SK H}\SK G$ is a monomorphism.
But, after identifying the latter with $\op_I\SK G$, this map is
just the inclusion $\SK H\subset\SK G$ on each coordinate and is
thus injective. Moreover this shows that the $\mathcal{K}G$--rank
of
 $M\ox_{\mathbb{Z} H}\SK G$ is also equal to the
cardinality of $I$.
\end{proof}

Finally, returning to the proof of our assertion, we claim that the map, $\phi\ox\id$, shown on the right hand side of the diagram above, is also injective. This will follow
immediately from Proposition~\ref{2-connected} below (setting $\G = B_n$), once we identify the
domain of $\phi\ox\id$ with $H_1(A;\SK B_n)$ and its range with $H_1(B;\SK B_n)$. The latter is immediate since we have previously observed that $\SK B_n$ is a flat $\mathbb{Z} B_n$--module. For the former, note that, since any module over a division ring is free, $\SK B_n$ is a free, and hence flat, $\SK A_n$ module. Moreover $\SK A_n$ is a flat $\mathbb{Z} A_n$--module. Hence
 $H_1(A;\SK B_n)\cong H_1(A;\mathbb{Z} A_n)\ox_{\mathbb{Z} A_n}\SK B_n$.  This completes the proof of the first claim of Theorem~\ref{main}, modulo the proof of
 Proposition~\ref{2-connected}.

The following proposition is an important result in its own right.
In the case that $\phi\co A\to B$ induces an isomorphism on $\hq1$,
this result was proved in \cite[Proposition 2.10]{COT}. The more general
result below seems to require a different proof. Here, by $\SK\G$ we mean
the right ring of fractions of the Ore domain $\bz\G$. Recall that any
homomorphism $\psi\co B\to\G$ endows $\bz\G$ and $\SK\G$ with left
$\bz B$--module structures. We call such a $\psi$ a \emph{coefficient
system} on $B$. If $\phi\co A\to B$ then
$\SK\G$ acquires a left $\bz A$--module structure via $\psi \circ \phi$.

\begin{prop}\label{2-connected} Suppose $\phi\co A\to B$ induces a
monomorphism (respectively, an isomorphism) on
$H_1(-;\mathbb{Q})$ and an epimorphism on $H_2(-;\mathbb{Q})$. Suppose also
that $A$ is finitely generated and $B$ is finitely related. Then for any coefficient system
$\psi\co B\to\G$, where $\G$ is a PTFA group, $\phi$ induces a
monomorphism (respectively, an isomorphism) $H_1(A;\SK\G)\to
H_1(B;\SK\G)$. Moreover, if the pair of Eilenberg--Maclane spaces $(K(B,1),K(A,1))$ has the homotopy type of a
relative $2$--complex, then the finiteness assumptions on $A$ and $B$ are not
necessary to get a monomorphism $H_1(A;\SK\G)\to H_1(B;\SK\G)$.
\end{prop}

Before proving Proposition~\ref{2-connected}, we finish the proof of the rest of the parts of our main theorem.

Note that the finiteness assumptions on $A$ and $B$ are used
only in the application of Proposition~\ref{2-connected}. Thus, if
the pair of Eilenberg--Maclane Spaces $(K(B,1),K(A,1))$ has the
homotopy type of a relative 2--complex, then we do not need these
finiteness assumptions to deduce the first part of the theorem.

Now, assume that $\phi$ induces an isomorphism on
$H_1(-;\mathbb{Q})$. We must show that
$\an/A\npp_{\sss H}\to\bn/B\npp_{\sss H}$ is a monomorphism between modules
of the same rank. The fact that this is a monomorphism follows
from the first part of the theorem. Since $\an/A\npp_{\sss H}$ and
$\an/[\an,\an]$ differ only by $\mathbb{Z}A_n$--torsion, they have
the same rank, $r_A$, as $\mathbb{Z}A_n$--modules. For the same
reason, $\an/A\npp_{\sss H}\ox_{\mathbb{Z}A_n}\SK(B_n)$ and
$\an/[\an,\an]\ox_{\mathbb{Z}A_n}\SK(B_n)$ are isomorphic. By
Lemma~\ref{mono}, the former has $\SK(B_n)$--rank equal to $r_A$
and hence so does the latter, which we have identified with
$H_1(A;\SK(B_n))$. If $\phi$ induces an isomorphism on
$H_1(-;\mathbb{Q})$ then note that $B$ must be finitely
generated. Hence Proposition~\ref{2-connected} applies to show
that $H_1(A;\SK(B_n))\cong H_1(B;\SK(B_n))$. Thus the latter has
$\SK(B_n)$--rank equal to $r_A$. But by applying the same reasoning
as above, we see that it has $\SK(B_n)$--rank equal to $r_B$, the
$\mathbb{Z}B_n$--rank of $\bn/B\npp_{\sss H}$.

If $\phi$ is also onto then it induces epimorphisms on all the
quotients by the terms of the torsion-free derived series. This
combined with the first part of the theorem implies that it
induces isomorphisms on all these quotients for $n\leq \omega$. Note also that since
this argument does not use the argument of the preceding
paragraph, it holds without the finiteness assumptions on $A$ and
$B$ if (B,A) is 2--dimensional in the sense described.

This concludes the proof of our main theorem, modulo the proof of Proposition~\ref{2-connected}.
\end{proof}

\begin{proof}[Proof of Proposition~\ref{2-connected}] We need the
following extension of a result of Ralph Strebel.

\begin{lem}\label{strebel} Suppose $\tl f\co M\to N$ is a homomorphism
between free $\bz\G$--modules with $\Gamma$  PTFA and let $f=\tl
f\ox\id$ be the induced homomorphism of abelian groups
$M\ox_{\bz\G}\bz\to N\ox_{\bz\G}\bz$. Then
$\rank_{\SK\G}(\image\tl f)\ge\rank_\mathbb{Q}(\image f)$.
\end{lem}

In \cite[page 305]{Str}, Strebel shows that, under the hypotheses of
Lemma~\ref{strebel}, if $f$ is injective then $\tilde{f}$ is
injective. There he shows that the class, $\mathcal{D}(\bz)$, of
groups $\Gamma$ for which this property is satisfied includes
torsion-free abelian groups and is closed under various natural
operations such as extensions. Consequently, any PTFA group is in
this class. This class of groups was previously called
$\textit{conservative}$ and was later shown by J Howie and H
Schneebli to coincide with the class of locally indicable groups
\cite{HoS}. Our lemma shows, in the case that $M$ is finitely generated,
 that more generally the rank of the kernel of $\tilde{f}$ is
at most the rank of the kernel of $f$.

\begin{proof}[Proof of Lemma~\ref{strebel}] By the rank of a homomorphism
we shall mean the rank of its image. Suppose that
$\rank_\mathbb{Q} f\geq r\le\infty$. Then there is a {\it
monomorphism} $g\co \bz^r\to N\ox_{\bz\G}\bz$ whose image is a
subgroup of $\image f$. If $e_i$, $1\le i\le r$ is a basis of
$\mathbb{Z}^r$, choose $M_i\in M\ox_{\mathbb{Z}\G}\mathbb{Z}$ such
that $f(M_i)=g(e_i)$. Since the ``augmentation'' $\e_M\co M\to
M\ox_{\mathbb{Z}\G}\mathbb{Z}$ is surjective there exist elements
$m_i\in M$ such that $\e_M(m_i)=M_i$. Consider the map $\tl
g\co (\mathbb{Z}\G)^r\to N$ defined by sending the $i^{\supth}$ basis
element to $\tl f(m_i)$. The augmentation of $\tl g$, $\tl
g\ox\id$, is the map
$(\mathbb{Z}\G)^r\ox_{\mathbb{Z}\G}\mathbb{Z}\to
N\ox_{\mathbb{Z}\G}\mathbb{Z}$ that sends $e_i$ to $\e_N(\tl f(m_i))=
f(\e_M(m_i))=g(e_i)$ and thus is seen to be
identifiable with $g$. In particular $\tl g\ox\id$ is a
monomorphism, and thus by \cite[page 305]{Str} (as mentioned above),
$\tl g$ is a monomorphism. Since the image of $\tl g$ lies in the
image of $\tl f$, $\tl g$ yields a monomorphism from
$(\mathbb{Z}\G)^r$ into the image of $\tl f$, showing that the
rank of $\image\tl f$ is at least $r$.
\end{proof}

Now we continue with the proof of Proposition~\ref{2-connected}.
We can find connected CW--complexes $X_A$, $X_B$ such that
$\pi_1(X_A)\cong A$ and $\pi_1(X_B)\cong B$ and whose universal
covers are contractible (classifying spaces for $A$, $B$). We can
find a cellular map $h\co X_A\to X_B$ inducing $\phi$ on $\pi_1$, and
by replacing $X_B$ by the mapping cylinder of $h$, we may assume
that $h$ embeds $X_A$ as a subspace of $X_B$. The finiteness
hypotheses on $A$ and $B$ are designed to ensure that (by proper
choice of the cell structure on $X_A$, $X_B$) we may assume that
the relative cellular chain group $C_2(X_B,X_A)$ (respectively,
$C_2(X_B,X_A)$ {\it and} $C_1(X_B,X_A)$) has finite rank. (Note that if $\phi$ induces an isomorphism
on $H_1(-;\mathbb{Q})$ then $B$ must be finitely presented). The
coefficient systems $\psi$ and $\psi\circ\phi$ induce covering
spaces $\wt X_B$ and $\wt X_A$ equipped with induced $\G$--actions
(these principal $\G$--bundles are connected covering spaces in case $\psi$ and
$\psi\circ\phi$ are surjective). We can lift the cell structure
and consider the relative free $\mathbb{Z}\G$--chain complex
$C_*(\wt X_B,\wt X_A)$ where $C_p(\wt X_B,\wt X_A)$ has
$\mathbb{Z}\G$--rank equal to the $\mathbb{Z}$--rank of
$C_p(X_B,X_A)$ and where the projection maps $\pi_p\co C_p(\wt
X_B,\wt X_A)\to C_p(X_B,X_A)$ commute with the $\p$ maps
$(\p_p\circ\pi_p=\pi_{p-1}\circ\tl\p_p)$. Moreover the projection
map can be identified with the augmentation $C_p(\wt X_B,\wt
X_A)\to C_p(\wt X_B,\wt X_A)\ox_{\mathbb{Z}\G}\mathbb{Z}\cong
C_p(X_B,X_A)$ in such a way that $\p_p=\tl\p_p\ox\id$. By Lemma
\ref{strebel} applied to $\tl\p_i\co C_i(\wt X_B,\wt X_A)\to
C_{i-1}(\wt X_B,\wt X_A)$ we have that
$\rank_{\SK\G}\tl\p_p\ge\rank_\mathbb{Q}\p_p$ for $p=1,2,3$. Since
the hypotheses imply that $H_2(X_B,X_A;\mathbb{Q})=0$
(respectively, in addition that $H_1(X_B,X_A;\mathbb{Q})=0$) we see  that $\rank_\mathbb{Q} H_2(C_*(X_B,X_A))=0$ (respectively, in
addition $\rank_\mathbb{Q} H_1(C_*(X_B,X_A))=0$). We now claim that
$\rank_{\SK\G}H_2(C_*(\wt X_B,\wt X_A))=0$ (respectively, in
addition $\rank_{\SK\G}H_1(C_*(\wt X_B,\wt X_A))$ is $0$). This follows
since (letting $c_2<\infty$ be the rank of $C_2(\wt X_B,\wt
X_A)$):
\begin{align*}
\rank_{\SK\G}(H_2(C_*(\wt X_B,\wt X_A))) &=
\rank_{\SK\G}(\ker\tl\p_2)
- \rank_{\SK\G}\tl\p_3\\
&= (c_2 - \rank_{\SK\G}\tl\p_2) - \rank_{\SK\G}\tl\p_3\\
&\le c_2 - \rank_\mathbb{Q}\p_2 - \rank_\mathbb{Q}\p_3\\
&= \rank_\mathbb{Q}(\ker\p_2) - \rank_\mathbb{Q}\p_3\\
&= \rank_\mathbb{Q}(H_2(C_*(X_B,X_A)))\\
&= 0.
\end{align*}
A similar argument holds for $H_1$ when $C_1(X_B,X_A)$ has finite
rank. Hence $H_2(\wt X_B,\wt X_A)$ has zero rank (respectively, in
addition $H_1(\wt X_B,\wt X_A)$ has zero rank). But the
equivariant homology modules $H_p(\wt X_B,\wt X_A)$ are well known
to be isomorphic to the homology with coefficients in $\psi$,
$H_p(X_B,X_A;\mathbb{Z}\G)$ \cite[Theorems VI3.4 and 3.4*]{Wh}.
Moreover the homology of a group $G$ is well known to be the same as
that of its associated Eilenberg--Maclane space $K(G,1)$ \cite[page 335]{HS}. Thus we have shown that
$H_2(B,A;\mathbb{Z}\G)\ox\SK\G=H_2(B,A;\SK\G)=0$ (respectively, in
addition $H_1(B,A;\SK\G)=0$) and consequently $\phi$ induces a
monomorphism (respectively an isomorphism) $\phi\co H_1(A;\SK\G)\to
H_1(B;\SK\G)$ as desired.

If we do not have the finiteness assumptions, but $(X_B,X_A)$ is a 2--complex, we can show that
$H_2(B,A;\mathbb{Q}\Gamma)=0$
by a direct application of Strebel's Lemma. Specifically, the hypothesis that $H_2(B,A;\mathbb{Q})=0$ implies that
the map $\partial _2:$\break$C_2(X_B,X_A;\mathbb{Q})\to C_1(X_B,X_A;\mathbb{Q})$ is injective. By Strebel's result, $\tilde{\partial
_2}$ is injective
and the claimed result follows. Hence, in this case even without the finiteness assumptions, we deduce that
$\phi\co H_1(A;\SK\G)\to H_1(B;\SK\G)$ is injective.
\end{proof}

\begin{cor}\label{free} Let  $F$ be a free group and $B$ be finitely related
with $H_2(B;\mathbb{Q})$ $=0$. Suppose $\phi\co F\to B$ induces
a monomorphism on $\hq1$. Then, for any finite $n$, it induces
monomorphisms $F/F\2n\subset B/B\2n$.
\end{cor}

\begin{proof}[Proof of Corollary~\ref{free}] Let $F$ be a finitely generated
free group.   Thus, by Theorem~\ref{main}, $F/F^{\sss (n)}\hookrightarrow
B/B^{\sss (n)}_{\sss H}$. Since $B^{\sss (n)}\subset B^{\sss (n)}_{\sss H}$, we see that
$F/F^{\sss (n)}\hookrightarrow B/B^{\sss (n)}$. Then we only need comment on
why it is not necessary to require that $F$ be finitely generated.
If it is not finitely generated, and if for some fixed $n$, the
map $F/F\2n\lra B/B\2n$ has a nontrivial element $w$ in its kernel, then $w$
is represented by a finite word which thus lies in a finite rank free
subgroup $G$ of $F$ (free on the letters that appear in $w$). Then $\phi$
induces a map $\phi\co G\to B$ that is a monomorphism on $\hq1$ and for which
the induced map the map $G/G\2n\lra B/B\2n$ has a nontrivial kernel,
contradicting Theorem~\ref{main}.
\end{proof}
\begin{rem}\label{noniso}

The first and second parts of Theorem~\ref{main} cannot be
improved to have isomorphisms in the conclusion, even for integral homology equivalences. For there exist
finitely presented groups $E$ and $\phi\co F\<x,y\>\to E$ inducing
isomorphisms on all integral homology groups but where the induced
map:
$$
\f{F^{\sss (1)}_{\sss H}}{F^{\sss (2)}_{\sss H}}\cong\f{F^{\sss (1)}}{F^{\sss (2)}}\lra
\f{E^{\sss (1)}}{E^{\sss (2)}}\cong\f{E^{\sss (1)}_{\sss H}}{E^{\sss (2)}_{\sss H}}
$$
is not surjective. One such example is $E=\<t,w,z\mid
t=z^3w^2tw^{-1}z^{-3}\>$ and $F=F\<t,z\>$. Then $F^{\sss (1)}/F^{\sss (2)}$
is a free $\bz[t^{\pm1},z^{\pm1}]$--module of rank 1 but
$E^{\sss (1)}_{\sss H}/E^{\sss (2)}_{\sss H}$ is not even a projective module. Such groups
arise commonly as the fundamental groups of the exteriors in $B^4$
of a set of ribbon disks for a ribbon link.
\end{rem}
\begin{rem}\label{nonfg}
Theorem~\ref{main} can
fail if $A$ is not finitely generated. Let $A=\langle x,w_i, i\in \mathbb{Z}\mid w_i=[x^{-1},w_{i+1}]\rangle$ and $B=\bz$. The abelianization $\phi\co A\to B$ induces an isomorphism on $H_1(-;\mathbb{Z})$ and an epimorphism on $H_2(-;\mathbb{Z})$, but $A^{\sss (1)}/A^{\sss (2)}$
has rank 1 (it has a $\mathbb{Z}[x^{\pm 1}]$--module presentation $\langle w_i, i\in \mathbb{Z} | w_i=w_{i+1}(x-1) \rangle$) so $\phi$ does not induce a monomorphism on
$A/A^{\sss (2)}_{\sss H}$.
\end{rem}
\begin{rem}\label{nonfg2}
The epimorphism part of the conclusion of
Proposition~\ref{2-connected} can fail if $B$ is not finitely
related using the same groups as above but with the roles
reversed, $\phi\co \bz\to A$. The same example shows that the part of
the conclusion of Theorem~\ref{main} about ``having the same
rank'' can fail if $B$ is not finitely related.
\end{rem}

\section{Homological completions and localizations}

In this section we construct a rational homological localization functor, $G\to \wt G$, that we call \emph{the torsion-free-solvable completion}. As we explain below, \emph{in the context of rational homological localization}, this can be viewed as an analogue of the Malcev completion wherein one replaces the lower central series by the torsion-free derived series. The latter is quite a bit more complicated because whereas the lower-central series quotients $G_n/G_{n+1}$ are \emph{trivial} modules (over $G/G_n$), their analogues $G^{\sss (n)}/G^{\sss (n+1)}$ are not. 
We parallel A. Bousfield's discussion of the Malcev completion (also called Malcev localization) \cite{B} and other homological localizations. At the end of the section we also compare our localization to the universal (integral) homological localization functor due to P Vogel and J Levine.

We remark that, in the context of \emph{linear algebraic groups}, C
Miller and R Hain have defined a ``relative solvable completion''
\cite{M, H}. This notion differs from ours. It is related to the lower-central series of the commutator subgroup. In particular, it is not invariant under homological equivalence.
We also note that, subsequent to our work, Jae Choon Cha has shown that our $\wt{G}$ is not the ``initial functor'' satisfying (1) and (2) below and has identified the initial such functor \cite{Cha}. 

\begin{thm}\label{universal} For any group $G$ there is a group $\wt G$ and
a homomorphism $f\co G\to\wt G$ such that:
\begin{enumerate}
\item $\ker f=G^{(w)}_{\sss H}$
\item If $A$ is finitely generated, $B$ is finitely-presented and $\phi\co A\to B$
induces an isomorphism (respectively, a monomorphism) on $\hq1$ and an epimorphism on $\hq2$ then
there is an isomorphism (respectively, monomorphism) $\tilde\phi\co \wt A\to\wt B$ such that the
following commutes:
\begin{equation*}
\begin{CD}
A          @>\phi>>   B\\
@Vf^AVV             @VVf^B V\\
\wt A         @>\tilde\phi>>   \wt B\\
\end{CD}
\end{equation*}
\end{enumerate}
\end{thm}

Before proceeding we review the Malcev completion so that the reader
can see the analogy to the torsion-free-solvable completion. Recall
that a nilpotent group $N$ is a \emph{uniquely divisible nilpotent
  group} if, for every positive integer $m$, the function $N$ to $N$
given by $x\to x^m$, is a bijection 
\cite[page 3]{B}. (If $\bigcap_{n=1}^{\infty}N^r_n=0$ then this can be
shown to be equivalent to requiring that the quotients,
$N^r_n/N^r_{n+1}$, of successive terms of the rational lower central
series are uniquely divisible $\mathbb{Z}$--modules). For a nilpotent
group, $N$, the \emph{Malcev completion} $N\to$MC($N$), usually
(sloppily) denoted $N \otimes \mathbb{Q}$, is a uniquely divisible
nilpotent group such that $N\to N \otimes \mathbb{Q}$ is initial among
maps to uniquely divisible nilpotent groups \cite[page 3]{B}
\cite{Q}\cite[Proposition 3.3]{PS}. This is identical to Bousfield's
$H\mathbb{Q}$-\emph{completion} \cite[Proposition 1.6]{B}. More
intuitively, $N\otimes \mathbb{Q}$ is obtained from the trivial group
by successive central extensions by the vector spaces
$(N^r_i/N^r_{i+1})\otimes \mathbb{Q}$ where the extensions are
``compatible'' with those that define $N$. For a general group , $G$,
the Malcev completion $MC(G)$ is $\varprojlim MC(G/G_n)$ or,
equivalently, $\varprojlim MC(G/G^r_n)$ (the
$\mathbb{Q}$-\emph{completion of G} in Bousfield's language \cite[page
  16]{B}).  Note that, by Stallings' Integral Theorem, the \emph{nilpotent completion}, $\varprojlim G/G_n$, is invariant under integral homology equivalence of groups but is \emph{not} invariant under rational homology equivalence. Even the \emph{rational nilpotent completion}, $\varprojlim G/G_n^r$, is not an invariant of rational homology equivalence (see Stallings' Rational Theorem). But the Malcev completion of $G$ \emph{is} invariant under rational homology equivalence (as are each of the groups $MC(G/G_n)$) as indicated by the last part of Stallings' Rational Theorem.

Theorem~\ref{universal} will be proved by directly constructing an \emph{$n$--solvable} version of $\wt G$ called $\wt G_n$ and then setting $\wt G=\varprojlim\wt G_n$ analogous to the construction of the Malcev/Bousfield completion. By way of further analogy, we will see that the torsion-free-solvable completion of a group $N$ is obtained from the trivial group successive extensions by the tensor product of the modules $N^{\sss (i)}_{\sss H}/N^{\sss (i+1)}_{\sss H}$ with appropriate skew fields. The torsion-free-solvable completion, $\wt G$, and the individual groups $\wt{(G/\gn)}$ will be seen to be invariant under rational homology equivalence. However, we have been unable to prove the result, which we expect is true, that our $\wt G_n$ is initial in the appropriate sense (among functors satisfying the properties of Theorem~\ref{superharvey}. Without this fact, the analogy to the Malcev completion is incomplete. The problem seems to be a failure of functoriality (since the torsion-free-derived series is not fully invariant). The authors expect this to be repaired by the modifications of Remark~\ref{Cohnremark}.

\begin{defn}\label{deftorsionfreesolv} A group $A$ is \textsl{$n$--torsion-free-solvable} if $A^{\sss (n)}_{\sss H}=0$ and is \textsl{torsion-free-solvable} if it is $n$--torsion-free-solvable for some integer $n$.
\end{defn}
Note that any torsion-free-solvable group $N$ is obtained from the trivial group by successive extensions by the torsion-free modules $N_{\sss H}^{\sss (i)}/N^{\sss (i+1)}_{\sss H}$. In particular, any such group is poly-(torsion-free-abelian).

\begin{defn}\label{torsionfreetower}(Compare \cite[Section 12]{B})\qua A collection of groups $A_n$, $ n\geq 0$, and group homomorphisms $f_n, \pi_n$ , $n\geq 0$, as below:
$$
\begin{array}{cccccc}
A     &\overset{f_n}{\lra}    &  A_n   &\overset{\pi_n}{\lra}
& A_{n-1}\\
\end{array}
$$
compatible in the sense that $f_{n-1}=\pi_n\circ f_n$,  is a \textsl{torsion-free-solvable tower for A} if, for each $n$, $A_n$ is $n$--torsion-free-solvable and the kernel of $\pi_n$ is contained in $(A_n)^{(n-1)}_{\sss H}$. A truncated tower, $A_n\to \ldots \to A_0$, is called a \textsl{tower of height $n$}.
\end{defn}

\begin{defn}\label{divisiblemodule} A  right module $M$ over an integral domain $R$ is a \textsl{(uniquely) divisible $R$--module} if, for each $m\in M$ and non-zero $r\in R$, there exists some (unique) $m'\in M$ such that $m=m'r$.
\end{defn}

\begin{defn}\label{divisibletfsolvablegp} A torsion-free-solvable group $A$ is a \textsl{(uniquely) divisible tor\-sion-free-solvable} group if, for each $n$, $A^{\sss (n)}_{\sss H}/A^{\sss (n+1)}_{\sss H}$ is a (uniquely) divisible $\mathbb{Z}[A/A^{\sss (n)}_{\sss H}]$--module.
\end{defn}

\begin{thm}\label{superharvey} For any group $G$ there exists a torsion-free-solvable tower, $\{\wt G_n, f_n\co G\to\wt G_n, \pi_n\co \wt G_n\to\wt G_{n-1}\}$ such that:

\begin{enumerate}

\item $\ker f_n=\gn$
\item  $\wt G_n$ is a uniquely divisible $n$--torsion-free-solvable group.

\item If $A$ is finitely generated,
$B$ is finitely presented and $\phi\co A\to B$ induces an isomorphism (respectively, monomorphism)
on $\hq1$ and an epimorphism on $\hq2$ then there is an
isomorphism (respectively, monomorphism) $\tilde\phi_n\co \wt A_n\to\wt B_n$ such that the following
commutes:
$$
\begin{array}{cccccc}
A     &\overset{f_n^A}{\lra}    &\wt A_n   &\overset{\pi_n^A}{\lra}
&\wt A_{n-1}\\
\text{\scriptsize$\phi$}\Big\downarrow         
&&\text{\scriptsize$\tilde\phi_n$}\Big\downarrow
&&\Big\downarrow\text{\scriptsize$\tilde\phi_{n-1}$}\\
B   &\overset{f^B_n}{\lra}   &\wt B_n   &\overset{\pi^B_n}{\lra}
&\wt B_{n-1} \\
\end{array}
$$
\item $\wt G_n$ depends only on $G/G^{\sss (n)}_{\sss H}$, that is if $\phi\co A\to B$ induces an isomorphism $A/A^{\sss (m)}_{\sss H}\cong B/B^{\sss (m)}_{\sss H}$ then it induces an isomorphism $\tilde \phi\co \wt A_n \to \wt B_n$ for all $n\leq m$. In particular, $\wt G_n \cong \wt {(G/G^{\sss (m)}_{\sss H})}_n$ if $n\leq m$.
\end{enumerate}
\end{thm}

\begin{proof}[Proof that Theorem~\ref{superharvey} implies
Theorem~\ref{universal}] Let $\wt G=\varprojlim\wt G_n$ and let
$f\co G\to\wt G$ be the map induced by the collection of $f_n\co G\to\wt
G_n$. Property $1$ of Theorem~\ref{universal} follows directly from property $1$ of Theorem~\ref{superharvey}.
\end{proof}

\begin{proof}[Proof of Theorem~\ref{superharvey}] We construct
such groups $(\wt G_n,f_n)$ recursively. Set $\wt G_0=\{ e\}$. Suppose
$\wt G_i$, $\pi_i$ and $f_i$ have been defined for $0\le i\le n$ in such a
way that Properties $1_n$ and $2_n$ are satisfied.

We first define $\wt G_{n+1}$.
Since $\wt G_{n}$ is torsion-free-solvable it is poly-torsion-free-abelian so $\bz\wt G_n$ is an Ore domain as in Proposition~\ref{PTFA}. Thus $\SK\wt G_n$
exists and is a $\bz\wt G_n-\SK\wt G_n$ bimodule. In particular $\SK\wt G_n$
is a $\bz G-\bz\wt G_n$ bimodule via $f_n\co G\to\wt G_n$ and so $H_1(G;\SK\wt
G_n)$ is defined and has the structure of a right $\bz\wt G_n$--module.
Let $\wt G_{n+1}$  be the semidirect product of $\wt G_n$ with
$H_1(G;\SK\wt G_n)$. Then we have the exact sequence
$$
1 \lra
H_1(G;\SK\wt G_n)\overset{i_n}{\lra}\wt
G_{n+1}\overset{\pi_{n+1}}{\lra}\wt G_n\lra 1.
$$

\begin{rem}\label{Bousfield} This construction is precisely analogous to Bousfield's construction of the $H\mathbb{Q}$--tower $\{f_n\co G\to T_n, n\geq 0\}$ for $G$, which goes as follows \cite[Section 3.4]{B}. Given $f_{n}\co G\to T_{n}$, let $V_{n}$ denote $H_2(f_{n};\mathbb{Q})$. Then the fundamental class in $H^2(f_{n};V_{n})$ determines a central extension
$$
1\lra V_{n}\lra T_{n+1}\lra T_{n}\lra 1.
$$
In our case, by analogy, given $f_{n}\co G\to \wt G_{n}$, let $K$
denote $H_2(f_{n};\SK\wt G_{n})$ where we take into account the module
structure. Now if we consider the exact sequence in homology for the
pair $(G,\wt G_{n})$ with coefficients in $\SK\wt G_{n}$, note that
$K\cong H_1(G;\SK\wt G_n)$ since $H_*(\wt G_{n};\SK\wt G_{n})=0$. Then
an extension of $\wt G_{n}$ by $K$ is determined as above. In fact we
shall see in Lemma~\ref{injective} that $K$ is an injective $\bz\wt
G_n$--module implying that $H^2(\wt G_n;K)=0$, and thus that the
semi-direct product is the unique such extension by $K$ \cite[page
  189, VI Theorem 10.3]{HS}.
\end{rem}

Next we want to show that $H_1(G;\SK\wt G_n)$ is a uniquely divisible $\bz\wt G_n$--module.
\begin{lem}\label{injectdivis} If $R$ is a right Ore Domain then a torsion-free right $R$--module $M$ is (uniquely) divisible if and only if it is a (uniquely) injective module (by \emph{(uniquely) injective} we mean that, for every monomorphism $\phi\co  L\to N$ and any homomorphism $\psi\co  L \to M$ there exists a (unique) extension $\psi '\co N\to M$.
\end{lem}
\begin{proof}[Proof of Lemma~\ref{injectdivis}] This follows from
  \cite[Proposition 3.7 page 58, Proposition  6.5 page 21]{Ste}. There Stenstr\"{o}m does not discuss the uniqueness condition, but this is easily verified by examining the proofs.
\end{proof}

\begin{lem}\label{injective} $H_1(G;\SK\wt G_n)$ is a uniquely injective (and uniquely divisible)
$\bz\wt G_n$--module, and is a uniquely injective (and uniquely divisible) $\bz[G/\gn
]$--module (via
$G/\gn\overset{f_n}{\lra}\wt G_n)$.
\end{lem}

\begin{proof}[Proof of Lemma~\ref{injective}] Since $H_1(G;\SK\wt G_n)$
is a right $\SK\wt G_n$--module and any $\SK\wt G_n$--module is free, it
suffices to know that $\SK\wt G_n$ itself is a uniquely injective
$\bz\wt G_n$--module \cite[Lemma 6.4 page 20]{Ste}). By Lemma~\ref{injectdivis} it suffices to see that $\SK\wt G_n$ is uniquely divisible. But this property of a classical quotient field is trivial to check (see
\cite[Proposition 3.7 page 58]{Ste}).

Since $f_n$ induces a monomorphism $\bz[G/\gn]\to\bz\wt G_n$ by $1_n$,
$\SK(G/\gn)$ embeds in $\SK\wt G_n$. Thus the $\bz[G/\gn]$
structure on $H_1(G;\SK\wt G_n)$ factors through $\SK(G/\gn)$.
Hence $H_1(G;\SK\wt G_n)$ is a $\SK(G/\gn)$--module, and therefore is
isomorphic to a direct sum of copies of $\SK(G/\gn)$. Since
$\SK(G/\gn)$ is itself a uniquely injective $\bz[G/\gn]$--module, as noted
above (since $\bz[G/\gn]$ is an Ore Domain), then $H_1(G;\SK\wt
G_n)$ is a uniquely injective $\bz[G/\gn]$--module.
\end{proof}

Now, application of the following Lemma ensures that $\wt G_{n+1}$ is a uniquely divisible $(n+1)$--torsion-free-solvable group, completing the verification of Property $2_{n+1}$.

\begin{lem}\label{semidirect} Suppose $B=A\rtimes C$ where $A$ is a
uniquely divisible $\bz C$--module and $C$ is $n$--torsion-free-solvable. Then $B$ is $(n+1)$--torsion-free-solvable and $A=B^{\sss (n)}_{\sss H}$. If $C$ is a uniquely divisible $n$--torsion-free-solvable group then $B$ is a uniquely divisible $(n+1)$--torsion-free-solvable group.
\end{lem}

\begin{proof}[Proof of Lemma~\ref{semidirect}] We may assume that $n$ is the least positive integer such that
$C^{\sss (n)}_{\sss H}=0$ since if for some $m<n$, $B^{\sss (m)}_{\sss H}=0$ then it follows that
$\bn=0$. We have the split exact sequence below
$$
1\lra A\overset{i}{\lra}B\overset{\pi}{\lra}C\lra 1.
$$
where $s\co C\to B$ is the splitting
map.

First we show by induction that $\pi_m\co B/B^{\sss (m)}_{\sss H}\to C/C^{\sss (m)}_{\sss H}$ exists and is an isomorphism for $0\le m\le
n$. This is trivially true for $m=0$. We assume it is true for all values at most some $m\le n-1$ and
establish it for $m+1$. Viewing $m$ as fixed, we first wish to show that
$s_m\co C/C^{\sss (m)}_{\sss H}\lra B/B^{\sss (m)}_{\sss H}$ exists and is an isomorphism. This requires a short induction. If $s_i$, $i<m$, exists and is an isomorphism
then by Proposition~\ref{functoriality}, $s_{i+1}$ exists. Since $i+1\le m$, our inductive hypothesis holds and so
$\pi_{i+1}$ exists and is an isomorphism. Moreover $\pi_{i+1}\circ s_{i+1}=\id$
so $s_{i+1}$ is an isomorphism. This completes the inductive proof that $s_m$ exists and is an
isomorphism. Returning to our proof that $\pi_{m+1}$ exists and is an isomorphism, choose a non-zero $c\in C^{(n-1)}_{\sss H}$. Since $m\le n-1$,
$c\in C^{\sss (m)}_{\sss H}$. Then $c-1$ is a non-zero element of $\bz C$ and since $A$
is a divisible $\bz C$--module, for any $a\in A$, there is an $\alpha\in A$ such
that $a=\alpha (c-1)$. When written in terms of $B$, this says that
$i(a)=[s(c)^{-1},i(\a)]$. Now note that since $i(A)=\ker\pi$ and $\pi_m$ is
injective, $i(A)\subset B^{\sss (m)}_{\sss H}$. Hence $i(\a)\in B^{\sss (m)}_{\sss H}$. Since $c\in
C^{\sss (m)}_{\sss H}$ and $s_m$ exists, $s(c)\in B^{\sss (m)}_{\sss H}$. Thus $i(a)\in
B^{\sss (m+1)}_{\sss H}$ and so $i(A)\subset B^{\sss (m+1)}_{\sss H}$. By Proposition
\ref{idempotency}, $\pi_{m+1}$ is an isomorphism. This concludes the
inductive proof that $\pi_{n}\co B/B^{\sss (n)}_{\sss H}\cong C/C^{\sss (n)}_{\sss H}$. However since
$C^{\sss (n)}_{\sss H}=0$, this implies that $i(A)=B^{\sss (n)}_{\sss H}$. Therefore $B^{\sss (n)}_{\sss H}$ is
abelian and $i(A)=B^{\sss (n)}_{\sss H}/[B^{\sss (n)}_{\sss H},B^{\sss (n)}_{\sss H}]$. This
quotient is naturally a $\bz[B/B^{\sss (n)}_{\sss H}]$--module (via conjugation as usual)
and the action factors through the isomorphism
$\pi_{n-1}\co \bz[B/B^{\sss (n)}_{\sss H}]\lra\bz[C/C^{\sss (n)}_{\sss H}]\cong\bz[C]$ (since $i(A)$ acts trivially). By hypothesis
$A$ is a uniquely divisible, hence torsion-free $\bz C$--module, so the submodule
$B^{\sss (n)}_{\sss H}/[B^{\sss (n)}_{\sss H},B^{\sss (n)}_{\sss H}]$ is a torsion-free
$\bz[B/B^{\sss (n)}_{\sss H}]$--module. Then, by definition, $B^{\sss (n+1)}_{\sss H}=0$.

Now suppose $C$ is a uniquely divisible $n$--torsion-free-solvable group. Since $A=B^{\sss (n)}_{\sss H}$, Proposition
\ref{idempotency} ensures that $\pi_m\co B^{(m-1)}_{\sss H}/B^{\sss (m)}_{\sss H}\to C^{(m-1)}_{\sss H}/C^{\sss (m)}_{\sss H}$ is an isomorphism for each $m\leq n$. Therefore $B^{(m-1)}_{\sss H}/B^{\sss (m)}_{\sss H}$ is a uniquely divisible $\mathbb{Z}[C/C^{(m-1)}_{\sss H}]$--module. Since $\bz[B/B^{(m-1)}_{\sss H}]\cong \bz[C/C^{(m-1)}_{\sss H}]$, $B^{(m-1)}_{\sss H}/B^{\sss (m)}_{\sss H}$ is a uniquely divisible $\bz[B/B^{(m-1)}_{\sss H}]$--module. It remains only to verify that $B^{\sss (n)}_{\sss H}/B^{\sss (n+1)}_{\sss H}$ ($=B^{\sss (n)}_{\sss H}=A$) is a uniquely divisible $\mathbb{Z}[B/B^{\sss (n)}_{\sss H}]$--module. Since $\mathbb{Z}[B/B^{\sss (n)}_{\sss H}]\cong \mathbb{Z}C$ this is true by hypothesis.
\end{proof}

Next we define $f_{n+1}$. Since $\wt G_{n+1}$ is a semidirect product we can specify $f_{n+1}$ uniquely by setting $f_{n+1}=(d_n,f_n)$ where $d_n\co G\to A$ is
a derivation, where we henceforth abbreviate $A=H_1(G;\SK\wt G_n)$, and $G$ acts (on
the right) through $G\overset{f_n}{\lra}\wt G_n$  \cite[VI Proposition
5.3]{HS}. Any such map $f_{n+1}$ will satisfy
$f_n=\pi_{n+1}\circ f_{n+1}$, so it only remains to specify $d_n$.
To define $d_n$, consider the exact sequence
$$
1 \lra\gn\lra G\lra G/\gn\lra 1
$$
and the induced 5--term exact sequence \cite[VI Theorem 8.1]{HS}
\begin{equation*}
\begin{CD}
0 \to \der(G/\gn,A)\lra\der(G,A)\overset{\pi}{\lra}\ho_{\bz[G/\gn]}(\gn/[\gn,\gn],A)\lra\\
H^2(G/\gn,A)\lra H^2(G,A).
\end{CD}
 \end{equation*}
We describe a canonical element of $\ho_{\bz[G/\gn]}(\gn/[\gn,\gn],A)$ , that we call $(f_n)_*$. Namely, consider the canonical projection $p\co H_1(G;\bz[G/\gn])\to B$ where $B=H_1(G;\bz[G/\gn])/\text{Torsion}$ and the canonical injection $i\co B\to B\ox\SK\wt G_n \cong A$ and set $(f_n)_*=i\circ p$ viewed as a map from $\gn/[\gn,\gn]$ to $A$. Since $A$ is an injective $\bz[G/\gn]$--module, $H^2(G/\gn;A)=0$ \cite[(2.4) page 189]{HS}, and so we can choose a derivation $d_n\in\der(G,A)$ such that $\pi (d_n)=(f_n)_*$.  Let $f_{n+1}\co G\to \wt G_{n+1}$ be the induced group homomorphism. We now show that, once $f_n$ is
fixed, $f_{n+1}$ is unique up to post-composition with an isomorphism.
Suppose $f'_{n+1}=(d'_n,f_n)$ is another homomorphism such that $\pi (d'_n)=(f_n)_*$. Then $d'_n=k + d_n$ where $k$ is the image of an element of $\der(G/G^{\sss (n)}_{\sss H},A)$. Since $A$ is an injective $\bz[G/\gn]$ (right)--module, $H^1(G/\gn,A)=0$ so $k$ is the image of a principal derivation, that is, $k(g)=\a (1-g^{-1})$ for some $\a\in A$. It is easy to verify
that the automorphism $\wt k$ of $\wt G_{n+1}$ given by conjugation by
$\a$, sends $(a,x)$ to $(\a (1-x^{-1})+a,x)$. Since $f_{n+1}(g)=(d_n(g),f_n(g))$, we see that $\wt k \circ f_{n+1}(g)$ is $(\a (1-f_n(g)^{-1})+ d_n(g),f_n(g))=(k(g)+ d_n(g),g)$ (since $G$ acts via $f_n$). Hence we have shown that $f_{n+1}'(g)= d'_n(g),f_n(g))=\wt k \circ f_{n+1}(g)$. Thus
$f_{n+1}$ is unique up to post-composition with an automorphism of $\wt G_{n+1}$. This completes the definition of $\{\wt G_{n+1},f_{n+1}\}$.

Having defined $f_{n+1}$, we must verify that $\{(\wt G_i,\wt f_i)| i\leq n+1\}$ is indeed a torsion-free-solvable tower of height $n+1$ for $G$. Since we have shown above that $\wt G_{n+1}$ is torsion-free-solvable, we only need that ker$(\pi_{n+1})\subset (\wt G_{n+1})^{\sss (n)}_{\sss H}$. But by Lemma~\ref{semidirect} these groups are equal.

We can now easily see that $f_{n+1}$
satisfies $1_{n+1}$, for if $f_{n+1}(g)=(0,e)$ then $g\in\ker f_n$
and $d_n(g)=0$. Thus $g\in \gn$, by  $1_{n}$, and $[g]$ lies in $\ker (f_n)_*$. But the kernel of
$(f_n)_*$ is precisely the subgroup represented by $G\npp_{\sss H}$, by
definition, so $g\in G\npp_{\sss H}$.

Now we verify Property $3$. Suppose $\phi\co A\to B$ is a homomorphism satisfying the
hypotheses of Property $3$ and  suppose, by induction
that $\phi$ induces an isomorphism (respectively, monomorphism) $\tilde\phi\co \wt A_n\to\wt B_n$
such that the diagram below commutes
\begin{equation*}
\begin{CD}
A          @>f_n>>   \wt A_n\\
@VV\phi V             @VV\tilde\phi_n V\\
B         @>f_n>>   \wt B_n\\
\end{CD}
\end{equation*}
Note that the composition:
$$
H_1(A;\SK\wt A_n)\overset{\tilde\phi_n}{\lra}H_1(A;\SK\wt
B_n)\overset{\phi_n}{\lra}H_1(B;\SK\wt B_n)
$$
is an isomorphism (respectively, monomorphism) of $\SK\wt A_n$ modules by
Proposition~\ref{2-connected}. This map, together with
$\tilde\phi_n$, induces an isomorphism (monomorphism) of extensions.
$$
\begin{array}{cccccc}
1     &\lra     &H_1(A;\SK\wt A_n)  &\lra\quad \wt A_{n+1}
&\overset{\pi_{n+1}}{\lra}\quad\wt A_n    &\lra 1\\
&&\quad\Big\downarrow\text{\scriptsize$(\tilde\phi_n)_*$}
&\qquad\ \ \Big\downarrow\text{\scriptsize$\tilde\phi_{n+1}$}
&\qquad \ \ \Big\downarrow\text{\scriptsize$\tilde\phi_n$}\\
1   &\lra   &H_1(B;\SK\wt B_n)  &\lra\quad\wt B_{n+1}
&\overset{\pi_{n+1}}{\lra}\quad\wt B_n    &\lra 1\\
\end{array}
$$
This finishes the proof of $3_{n+1}$ except for verifying the
compatibility of maps. To verify that $\tilde\phi_{n+1}\circ
f^A_{n+1}=f^B_{n+1}\circ\phi$, observe that since $\wt A_{n+1}$ and
$\wt B_{n+1}$ are semidirect products and since $\tilde\phi_n\circ
f^A_n=f^B_n\circ\phi$, it suffices to check $\tilde\phi_{n+1}\circ
f^A_{n+1}=f^B_{n+1}\circ\phi$ as maps from $\ker f^A_n=\an$ to
$H_1(B;\SK\wt B_n)$. For this recall that $f^A_{n+1}$ restricted to
$\ker f^A_n$ (and $f^B_{n+1}$ restricted to $\phi(\an)\subset\bn=\ker
f^B_n$) is given by the canonical map $(f^A_n)_*$ induced by $\phi$
and $f^A_n$. On the other hand, $\tilde\phi_{n+1}$ restricted to
$f^A_{n+1}(\an)\sbq H_1(A;\SK\wt A_n)$ is determined by $\phi$.
Details are left to the reader.

Finally, to establish Property $4$, suppose that $\phi\co A\to B$ induces an isomorphism $A/A^{\sss (m)}_{\sss H}\cong B/B^{\sss (m)}_{\sss H}$. Inductively we may assume that this isomorphism extends to isomorphisms $\wt A_n\cong \wt B_n$ and try to show that they extend to an isomorphism $\wt A_{n+1}\cong \wt B_{n+1}$ as long as $n+1\leq m$. By the Five Lemma it suffices to show that the modules $H_1(A;\SK\wt A_n)$ and $H_1(B;\SK\wt B_n)$ (modules over isomorphic rings) are isomorphic. But $H_1(A;\SK\wt A_n)$ is completely determined by $H_1(A;\mathbb{Z}[A/\an])/\{\Z[A/A^{\sss (n)}_{\sss H}]\text{--Torsion}\}$ together with the inclusion $A/\an \to \wt A_n$. But the former is precisely $\an/A^{\sss (n+1)}_{\sss H}$. If $n+1\leq m$ then all of these are carried isomorphically onto the corresponding modules for $B$.

This completes the proof of Theorem~\ref{superharvey}.
\end{proof}

\begin{cor}\label{stabilize} The torsion-free derived series stabilizes at $n$ if and only if $\{ \wt G_n\}$ stabilizes at $n$. In particular, if  $\gn=G^{\sss (n+1)}_{\sss H}$, then $\wt G_n=\wt G_{n+k}$ for all $k\geq 0$ and
$\wt G=\wt G_n\cong \wt{(G/\gn)}$. Also, for any group $G$,
$\wt{(G/\gn)}_n=\wt {(G/\gn)}$.
\end{cor}

\begin{proof}[Proof of Corollary~\ref{stabilize}] Suppose that $\gn=G^{\sss (n+1)}_{\sss H}$. Then by construction $\wt G_{n+1}= \wt G_n$. Since $\wt G_n$ stabilizes at $n$, by definition $\wt G=\wt G_n$. Conversely, if $\wt G_{n+1}= \wt G_n$ then the kernel of $G\to \wt G_{n+1}$ equals the kernel of $G\to \wt G_{n}$, thus by Theorem~\ref{superharvey}, $G^{\sss (n+1)}_{\sss
H}=\gn$. For the final claim, merely note that $G/\gn$ is
$n$--torsion-free-solvable by Proposition~\ref{idempotency}.
\end{proof}

It follows that the torsion-free-solvable completion may be constructed
analogously to the Malcev completion.

\begin{cor}\label{invlimit} The torsion-free-solvable completion $\wt G$ of $G$ is the inverse limit $\varprojlim \wt{(G/\gn)}$ of the torsion-free-solvable completions of the corresponding quotients by the torsion-free derived series.
\end{cor}

\begin{proof}[Proof of Corollary~\ref{invlimit}] Since $\wt G$ has been defined as $\varprojlim\wt G_n$, we must show that $\wt G_n=\wt{(G/\gn)}$. By Proposition~\ref{idempotency}, the map $G\to G/\gn$ induces an isomorphism $G/\gn \cong (G/\gn)/(G/\gn)_{\sss H}^{\sss (n)}$. Thus, by Theorem~\ref{superharvey} $(4)$, $\wt G_n=\wt{(G/\gn)}_n$ which in turn equals $\wt{(G/\gn)}$ by the last statement of Corollary~\ref{stabilize}.
\end{proof}

\begin{ex}\label{finite2} If $\beta _1(G)=0$, for example if $G$ is a finite group, then the torsion-free-solvable completion and the Malcev completion agree and are both trivial. For we observed in Example~\ref{finite} that the torsion-free derived series stabilizes at $n=0$ in such cases. Thus by Corollary~\ref{stabilize}, $\wt G\cong\wt G_0\cong 0\cong MC(G)$.
\end{ex}
\begin{ex}\label{abelian2} If $G$ is a  non-trivial torsion-free abelian group (equivalently $G^{\sss (1)}_{\sss H}=\{e\}$), then $\wt G$ again agrees with the Malcev completion. For $G^{\sss (1)}_{\sss H}=0$ so the torsion-free derived series stabilizes at $n=1$. Thus by Corollary~\ref{stabilize}, $\wt G\cong \wt G_1= G\otimes\mathbb{Q}\cong MC(G)$. Note that the algebraic closure of $\mathbb{Z}^m$ is $\mathbb{Z}^m$ but the Malcev completion and the torsion-free-solvable completion are $\mathbb{Q}^m$.
\end{ex}
\begin{ex}\label{nilpotent2} If $G$ is nilpotent then whereas $G/\text{torsion}$ embeds in MC($G$), the torsion-free-solvable completion of $G$ is just the abelianization of MC($G$). For we saw in Example~\ref{nilpotent} that in this case the torsion-free derived series stabilizes at $n\leq 1$. Thus, by Corollary~\ref{stabilize}, $\wt G=\wt G_1=\wt {(G/G^{\sss (1)}_{\sss H})}$ and the latter equals $G/G^{\sss (1)} \otimes \mathbb{Q}$ by Example~\ref{abelian2}. This is a flaw in the torsion-free derived series that would be corrected by Remark~\ref{Cohnremark}.
\end{ex}

\begin{ex}\label{bettione2} If $\beta _1(G)=1$ then the epimorphism $G\to \mathbb{Z}$ is rationally 2--connected and so, by Theorem~\ref{universal}, it induces isomorphisms $\wt G \cong \wt{\mathbb{Z}} \cong \mathbb{Q}\cong MC(G)$ by Example~\ref{abelian2}. More generally, if $G$ is a group whose classical Alexander module $G^{\sss (1)}_{\sss H}/[G^{\sss (1)}_{\sss H},G^{\sss (1)}_{\sss H}]$ has rank zero then, as noted in Example~\ref{rankzero}, the torsion-free derived series stabilizes at $n\leq 1$. Thus, by Corollary~\ref{stabilize}, $\wt G=\wt G_1=\wt {(G/G^{\sss (1)}_{\sss H})}=G/G^{\sss (1)} \otimes \mathbb{Q}$.
\end{ex}
\begin{ex}\label{freesolvable2}  If the rank of one of the modules $\gn /[\gn ,\gn ]$, is zero then the torsion-free derived series stabilizes at $n$ or less. Thus $\wt G=\wt G_n$ is solvable by Corollary~\ref{stabilize}. For example if $G$ is the free solvable group $F/F^{\sss (n)}$, then $G$ embeds in its torsion-free-solvable completion since $(F/F^{\sss (n)})^{\sss (\omega)}_{\sss H}=0$.
\end{ex}
\begin{ex}\label{general2} If $F$ is a non-abelian free group then, since the torsion-free derived series stabilizes at $\omega$, $\wt F$ is not nilpotent. In fact, since $F^{\sss (\omega)}_{\sss H}=0$, $F$ embeds in $\wt F$. But $\wt F$ is much larger, for if $G=\pi_1(S^3-L)$, where $L$ is a boundary link, then the meridional map $F\to G$ has a right inverse $G\to F$ as we saw in Example~\ref{general} and it follows from Theorem~\ref{universal} that $\wt G \cong \wt F$ and $G/G^{\sss (\omega)}_{\sss H}$ embeds in $\wt F$. Indeed if $G$ is any finitely-presented group with $H_1(G;\mathbb{Q})\cong \mathbb{Q}^m$ and $H_2(G;\mathbb{Q})=0$ then there is a map $F\to G$ that is 2--connected on rational homology and so, by Theorem~\ref{universal}, $G/G^{\sss (\omega)}_{\sss H}$ embeds in $\wt F$. In particular, as we showed in Theorem~\ref{algclosure}, $\wh F/(\wh F)^{\sss (\omega)}_{\sss H}$ embeds in $\wt F$.
\end{ex}

\subsection{Relations with the Levine--Vogel algebraic closure}

We now recall the (integral) homological localization functor of
Levine--Vogel and discuss the relations between it and the
torsion-free-solvable completion. J Levine defined the
\emph{algebraic closure}, $\hat f\co G\to \wh {G}$ of a group $G$,
unique up to isomorphism \cite[page 574]{L1}. He pointed out that, if $G$ is finitely presented, it coincides with a notion previously investigated by P. Vogel (see \cite{LD}). Specifically, if $X$ is a finite CW--complex with $G=\pi_1(X)$ then $\wh G$ is $\pi_1(EX)$ where $EX$ is the \emph{Vogel localization} of $X$. The Levine--Vogel algebraic closure is a universal (integral) homological localization in the following sense. If $A$ is finitely generated, $B$ is finitely presented and $\phi\co  A\to B$ is a homomorphism that induces an isomorphism on $H_1(-;\mathbb{Z})$, an epimorphism on $H_2(-;\mathbb{Z})$, and such that $B$ is the normal closure of $\phi(A)$, then there is an isomorphism $\hat\phi\co \wh A\to\wh B$. Moreover the algebraic closure is the ``initial'' group with this property. Since $\wt G$ also satisfies this property, we have a canonical map $\wh{G} \xrightarrow{\theta} \wt{G}$ that provides a factorization of $f$ (above) as $G \xrightarrow{\hat f} \wh{G} \xrightarrow{\theta} \wt{G}$. This allows us to relate $\wh G$ to $\wt G$.

\begin{thm}\label{algclosure}Suppose G is a finitely presented group. Then there is a canonical map from the algebraic closure of $G$ to the torsion-free-solvable completion of $G$, $\wh{G} \xrightarrow{\theta} \wt{G}$ whose kernel is $(\wh{G})_{\sss H}^{\sss (\omega)}$.
\end{thm}
\begin{proof} [Proof of Theorem~\ref{algclosure}] Recall from \cite[Proposition 6]{L1} that $\wh{G}$ can be expressed as $\underrightarrow{\lim} P_i$ where 
$G=P_0\to P_1\to \dots \to P_k\to \dots \to \wh G$, each $P_i$ is finitely presented and each map satisfies the Levine--Vogel conditions above. 
Suppose $x\in\ker\theta$. Then $x$ is the image of $p\in P_k$ for some $k$, and $p$ lies in the kernel of the induced map $P_k\to \wh G \to \wt G$. Since $\wt G \cong \wt{P_k}$ by Theorem~\ref{universal}, 
one concludes that $p$ is in the kernel of the canonical map $P_k \to \wt {P_k}$. Again by Theorem~\ref{universal}, $p\in (P_k)_{\sss H}^{\sss (\omega)}$. To now conclude that $x\in (\wh{G})_{\sss H}^{\sss (\omega)}$, finishing the verification that $\ker\theta\subset(\wh{G})_{\sss H}^{\sss (\omega)}$, we need the following.

\begin{lem}\label{directlimit} If $B= \underrightarrow{\lim} P_i$ where
$P_i$ are finitely presented and each map $\phi _i\co P_i\to P_{i+1}$ induces
a monomorphism on $H_1(-;\mathbb{Q})$ and an epimorphism on
$H_2(-;\mathbb{Q})$, then for each $i$ and each $n\leq\omega$, the
induced map $(P_i)/(P_i)^{\sss (n)}_{\sss H}\lra B/B^{\sss (n)}_{\sss H}$ is injective.
\end{lem}
\begin{proof} [Proof of Lemma~\ref{directlimit}] The proof is by induction
on $n$. Note that it is trivially true for $n=0$. The proof follows
exactly the proof of Theorem~\ref{main}. The Lemma would follow directly
if $B$ were finitely related. However this hypothesis on $B$ is only
needed in the inductive step to establish the monomorphism conclusion of
Proposition~\ref{2-connected} for the map $H_1(P_i;\SK \G) \to H_1(B;\SK
\G
)$ where $\G= B/\bn$. Since each $P_i$ is finitely presented,
Proposition~\ref{2-connected} shows that each map $H_1(P_i;\SK \G)\to
H_1(P_{i+k};\SK \G)$ is injective. Thus $H_1(P_i;\SK \G)\to
\underrightarrow{\lim} H_1(P_{i+k};\SK \G)$ is injective ($i$ fixed, $k$
varying). Since $B= \underrightarrow{\lim} P_{i+k}$ and since homology
commutes with direct limits \cite[page 121]{Br}, the conclusion of
Proposition~\ref{2-connected} holds for each map $P_i\to B$.
\end{proof}

Continuing with the proof of Theorem~\ref{algclosure},
now suppose $x\in (\wh{G})_{\sss H}^{\sss (\omega)}$. Then $x$ is the image of $p\in P_i$ for some $i$ and so, by the case $n=\omega$ of Lemma~\ref{directlimit}, $p\in (P_i)_{\sss H}^{\sss (\omega)}$. Hence, by Theorem~\ref{universal}, $p$ is in the kernel of the canonical map $P_i\to \wt{(P_i)}$ and since $\wt{(P_i)}\cong \wt G$ we have $x\in \ker\theta$ as desired.
 \end{proof}

These observations suggest that it might be profitable to consider the functor $G\to \wh G/(\wh{G})_{\sss H}^{\sss (\omega)}$ which can be described as a direct limit $\underrightarrow{\lim} P_i/(P_i)_{\sss H}^{\sss (\omega)}$. We say that a homomorphism is \textsl{2--connected on integral homology} or \textsl{2--connected} if it induces an isomorphism on $H_1(-;\mathbb{Z})$ and an epimorphism on $H_2(-;\mathbb{Z})$. We say that a homomorphism is \textsl{2--connected on rational homology} or \textsl{rationally 2--connected} if it induces an isomorphism on $H_1(-;\mathbb{Q})$ and an epimorphism on $H_2(-;\mathbb{Q})$. The interesting question is: If $A\to B$ is rationally 2--connected but fails to be 2--connected with integer coefficients or fails to have the Levine--Vogel normal generation condition, what can be said about the relationship between $\wh A$ and $\wh B$? We have the following result, used in \cite{Ha2} to prove that certain invariants of homology cobordism and link concordance are actually invariants of rational homology cobordism.

\begin{cor}\label{derivedlocal} If $A$ and $B$ are finitely-presented and $\phi\co A\to B$
induces an isomorphism on $\hq1$ and an epimorphism on $\hq2$ then $\phi$ induces an embedding $\wh A/(\wh A)_{\sss H}^{\sss (\omega)} \subset \wh B/(\wh B)_{\sss H}^{\sss (\omega)}$.
\end{cor}

\begin{proof}[Proof of Corollary~\ref{derivedlocal}] Since any algebraic closure $A\to \wh A$ is 2--connected with
integral coefficients \cite[Proposition 4]{L1}, $\phi$ induces an embedding $\wh A/(\wh A)_{\sss H}^{\sss (1)}$ $\subset \wh B/(\wh B)_{\sss H}^{\sss (1)}$, 
since these quotients are merely $H_1(-;\mathbb{Z})$ modulo torsion. We proceed by induction. Suppose that $\phi$ induces a monomorphism 
$\wh A/(\wh A)_{\sss H}^{\sss (n-1)} \subset \wh B/(\wh B)_{\sss H}^{\sss (n-1)}$. By Proposition~\ref{functoriality}, $\phi$ induces a homomorphism 
$\wh A/(\wh A)_{\sss H}^{\sss (n)} \to \wh B/(\wh B)_{\sss H}^{\sss (n)}$. It suffices to show that this is injective. Suppose $a\in \wh A$ such that 
$\hat \phi (a)=b\in (\wh B)^{\sss (n)}_{\sss H}$. As above, we know that $\hat{B}$ can be expressed as $\underrightarrow{\lim} Q_i$ where 
$B=Q_0\to Q_1\to \dots \to Q_k\to \dots \to \wh B$, each $Q_i$ is finitely presented and each map is 2--connected with integral coefficients. Similarly, $\wh A=\underrightarrow{\lim} P_i$. Suppose $q\in Q_i$ has image $b$ and $p\in P_j$ has image $a$. For simplicity we abbreviate $Q=Q_i$ and $P=P_j$. By Lemma~\ref{directlimit}, $q\in Q^{\sss (n)}_{\sss H}$. However, since $B\to Q$ and $A\to P$ are 2--connected and $A\to B$ is rationally 2--connected, $\wt B_n \cong \wt Q_n$, $\wt A_n \cong \wt P_n$ and $\wt A_n \cong \wt B_n$ by Theorem~\ref{superharvey}. Also by part $1)$ of this theorem, the image of $q$ under $Q\to \wt Q_n$ is zero and hence the image of $b$ under $\wh B\to \wt B\to \wt B_n$ is zero. It follows that the image of $a$ under the map $\wh A\to \wt A_n$ is zero. Consequently the image of $p$ under the map $P\to \wt P_n$ is zero and hence, by part $1)$ of Theorem~\ref{superharvey}, $p\in P^{\sss (n)}_{\sss H}$. Thus, by Lemma~\ref{directlimit}, $a\in (\wh A)^{\sss (n)}_{\sss H}$. Thus we have shown that $\wh A/(\wh A)_{\sss H}^{\sss (n)} \to \wh B/(\wh B)_{\sss H}^{\sss (n)}$ is injective as desired.
\end{proof}

\begin{rem}\label{Cohnremark} We are aware that a smaller series (that looks much less natural from an algebraic standpoint) than the torsion-free derived series seems to be
 \emph{more} natural from the point of view of rational homology equivalence. Note that $\mathbb{Z}[G/G^{\sss (1)}_{\sss H}]$ is a Laurent polynomial ring
  $\mathbb{Z}[x_1^{\pm 1},...,x_m^{\pm 1}]$. We could alternatively define $G^{\sss (2)}_{\sss H}$ to be the inverse image of not the \emph{full} torsion submodule 
  of $G^{\sss (1)}_{\sss H}/[G^{\sss (1)}_{\sss H},G^{\sss (1)}_{\sss H}]$, but rather the submodule of elements annihilated by some element of the set $S$ of polynomials 
  whose image under the augmentation map $\mathbb{Z}[x_1^{\pm 1},...,x_m^{\pm 1}]\to \mathbb{Z}$ is non-zero. Then, leaving the definition of the higher
   $G^{\sss (n)}_{\sss H}$ the same, leads to a series whose terms are smaller than the torsion-free derived series, still containing the derived series, and for which our main 
   theorem remains true because a more precise version of Proposition~\ref{2-connected} is known to hold (with $\G = \mathbb{Z}^m$ and) with the quotient field being replaced
    by the subring $S^{-1}\mathbb{Z}[G/G^{\sss (1)}_{\sss H}]$. Generalizing this, a ``better'' series, $G^{\sss (n)}_*$, can be defined as follows. Once $G^{\sss (n)}_*$ has 
    been defined, let $CL(G/G^{\sss (n)}_*)$ stand for the Cohn localization of the augmentation $\mathbb{Z}[G/G^{\sss (n)}_*]\to \mathbb{Q}$. By this we mean the initial ring 
    map $\mathbb{Z}[G/G^{\sss (n)}_*]\to CL(G/G^{\sss (n)}_*)$, such that any square matrix over $\mathbb{Z}[G/G^{\sss (n)}_*]$ whose augmentation is invertible over
     $\mathbb{Q}$ becomes invertible over $CL(G/G^{\sss (n)}_*)$. Note that in the case $n=1$, $CL(G/G^{\sss (1)}_{\sss *})$ is known to be just 
     $S^{-1}\mathbb{Z}[G/G^{\sss (1)}_{\sss *}]$ with $S$ as above, because since $\mathbb{Z}[G/G^{\sss (1)}_{\sss *}]$ is commutative, a matrix is invertible if and only if its
      determinant is invertible. We could then define $G^{\sss (n+1)}_*$ to be the kernel of the composition below:
$$
G^{\sss (n)}_* \xrightarrow{\pi_n} \f{G^{\sss (n)}_*}{[G^{(n)}_*,G^{(n)}_*]}=H_1(G;\mathbb{Z}[G/G^{\sss (n)}_*])\lra H_1(G;CL(G/G^{\sss (n)}_*)).
$$
This is in contrast to the torsion-free derived series, wherein $G^{\sss (n+1)}_{\sss H}$ is the kernel of the composition:
$$
G^{\sss (n)}_{\sss H} \xrightarrow{\pi_n} \f{\gn}{[\gn,\gn]}=H_1(G;\mathbb{Z}[G/\gn])\lra H_1(G;\SK(G/\gn)).
$$
The map from $\mathbb{Z}[G/G^{\sss (n)}_*]$ to its quotient field $\SK(G/G^{\sss (n)}_*)$ is known to factor through $CL(G/G^{\sss (n)}_*)$.
This series is then fully invariant. Much of the structure of our main theorem holds, but there are several problems involving flatness. The Cohn localization is not usually a flat module
 over $\mathbb{Z}[G/G^{\sss (n)}_*]$ and so we have not been able to duplicate the full strength of our main theorem. However we can prove that
$$
G^{(n)}_r\subset G^{(n)}_*\subset G^r_{2^n}
$$
and that if $\phi\co A\to B$ is an \emph{epimorphism} of groups which is rationally 2--connected then it induces isomorphisms $A/A^{\sss (n)}_*\cong B/B^{\sss (n)}_*$ for all $n$. The proofs of these results will appear in a subsequent paper.
\end{rem}


\begin{thebibliography}

\bibitem{B}
\textbf{K Bousfield}, \emph{Homological localization towers for groups and
  $\pi$--modules}, Mem. Amer. Math. Soc. 186 (1977) \MR{0447375}

\bibitem{Br}
\textbf{K\,S Brown}, \emph{Cohomology of groups}, Graduate Texts in Mathematics
  87, Springer--Verlag, New York (1982) \MR{672956}

\bibitem{Ca}\textbf{A\,J Casson}, \emph{Link Cobordism and Milnor's
Invariant},  Bull. London Math.  Soc. 7 (1975) 39--40
\MR{0362286}

\bibitem{Cha}
\textbf{JC Cha}, \emph{Injectivity Theorems and Algebraic Closure of Groups
  with Coefficients} Preprint

\bibitem{C}
\textbf{T\,D Cochran},
\href{http://www.maths.warwick.ac.uk/agt/AGTVol4/agt-4-19.abs.html}%
{\emph{Noncommutative knot theory}}, Algebr. Geom. Topol.
  4 (2004) 347--398 \MR{2077670}

\bibitem{COT}
\textbf{T\,D Cochran}, \textbf{K\,E Orr}, \textbf{P Teichner}, \emph{Knot
  concordance, {W}hitney towers and {$L\sp 2$}--signatures}, Ann. of Math. (2)
  157 (2003) 433--519 \MR{1973052}

\bibitem{COT2}
\textbf{T\,D Cochran}, \textbf{K\,E Orr}, \textbf{P Teichner}, \emph{Structure
  in the classical knot concordance group}, Comment. Math. Helv. 79 (2004)
  105--123 \MR{2031301}

\bibitem{CT}
\textbf{T\,D Cochran}, \textbf{P Teichner}, \emph{Knot Concordance and von Neumann
  $\rho$--invariants}, \arxiv{math.GT/0411057}

\bibitem{Co}
\textbf{P\,M Cohn}, \emph{Free rings and their relations}, London Mathematical
  Society Monographs 19, Academic Press Inc. [Harcourt Brace Jovanovich
  Publishers], London (1985) \MR{800091}

\bibitem{H}
\textbf{R\,M Hain}, \emph{Completions of mapping class groups and the cycle
  {$C-C\sp -$}}, from: ``Mapping class groups and moduli spaces of Riemann
  surfaces (G\"ottingen, 1991/Seattle, WA, 1991)'', Contemp. Math. 150, Amer.
  Math. Soc. Providence, RI (1993)  75--105 \MR{1234261}

\bibitem{Ha2}
\textbf{S Harvey}, \emph{Homology Cobordism Invariants of 3--Manifolds and the
  Cochran--Orr--Teichner Filtration of the Link Concordance Group}, preprint

\bibitem{Ha1}
\textbf{S\,L Harvey}, \emph{Higher-order polynomial invariants of 3--manifolds
  giving lower bounds for the {T}hurston norm}, Topology 44 (2005) 895--945
  \MR{MR2153977}

\bibitem{HS}
\textbf{P\,J Hilton}, \textbf{U Stammbach}, \emph{A course in homological
  algebra}, Graduate Texts in Mathematics 4, Springer--Verlag, New York (1997)
  \MR{1438546}

\bibitem{HoS}
\textbf{J Howie}, \textbf{H\,R Schneebeli}, \emph{Homological and topological
  properties of locally indicable groups}, Manuscripta Math. 44 (1983) 71--93
  \MR{709846}

\bibitem{Ki1}
\textbf{T Kim}, \emph{Filtration of the classical knot concordance group and
  {C}asson--{G}ordon invariants}, Math. Proc. Cambridge Philos. Soc. 137 (2004)
  293--306 \MR{2092061}

\bibitem{Ki2}
\textbf{T Kim}, \emph{An infinite family of non-concordant knots having the
  same {S}eifert form}, Comment. Math. Helv. 80 (2005) 147--155 \MR{MR2130571}

\bibitem{LD}
\textbf{J LeDimet}, \emph{Cobordism d'Enlacements de Disques}, Bull. Math. Soc.
  France No.32 116 (1988)

\bibitem{L1}
\textbf{J\,P Levine}, \emph{Link concordance and algebraic closure. {II}},
  Invent. Math. 96 (1989) 571--592 \MR{996555}

\bibitem{Le}
\textbf{J Lewin}, \emph{A note on zero divisors in group-rings}, Proc. Amer.
  Math. Soc. 31 (1972) 357--359 \MR{0292957}

\bibitem{M}
\textbf{C Miller}, \emph{Exponential iterated integrals and the relative
  solvable completion of the fundamental group of a manifold}, Topology 44
  (2005) 351--373 \MR{2114712}

\bibitem{PS}
\textbf{S Papadima}, \textbf{A\,I Suciu}, \emph{Chen {L}ie algebras}, Int.
  Math. Res. Not.  (2004) 1057--1086 \MR{2037049}

\bibitem{P}
\textbf{D\,S Passman}, \emph{The algebraic structure of group rings},
  Wiley-Interscience [John Wiley \& Sons], New York (1977) \MR{470211}

\bibitem{Q}
\textbf{D Quillen}, \emph{Rational homotopy theory}, Ann. of Math. (2) 90
  (1969) 205--295 \MR{0258031}

\bibitem{Ro}
\textbf{SK Roushon}, \emph{Topology of 3--manifolds and a Class of Groups II},
  Boletin de la Sociedad Matematica Mexicana 3a, Serie 10(3) (2004)

\bibitem{St}
\textbf{J Stallings}, \emph{Homology and central series of groups}, J. Algebra
  2 (1965) 170--181 \MR{0175956}

\bibitem{Ste}
\textbf{B Stenstr{\"o}m}, \emph{Rings of quotients}, Springer--Verlag, New York
  (1975) \MR{0389953}

\bibitem{Str}
\textbf{R Strebel}, \emph{Homological methods applied to the derived series of
  groups}, Comment. Math. Helv. 49 (1974) 302--332 \MR{0354896}

\bibitem{T}
\textbf{P Teichner}, \emph{Knots, von {N}eumann signatures, and grope
  cobordism}, from: ``Proceedings of the International Congress of
  Mathematicians, Vol. II (Beijing, 2002)'', Higher Ed. Press, Beijing (2002)
  437--446 \MR{1957054}

\bibitem{Wh}
\textbf{G\,W Whitehead}, \emph{Elements of homotopy theory}, Graduate Texts in
  Mathematics 61, Springer--Verlag, New York (1978) \MR{516508}

\end{thebibliography}
\end{document}